\documentclass[leqno]{amsart}
\usepackage{egastyle}
\usepackage{mymacros}
\usepackage[all,tips]{xy}
\usepackage{url}
\usepackage{bm}
\usepackage{verbatim}
\CompileMatrices

\def\cU{{\mathcal U}}
\def\cV{{\mathcal V}}

\def\B{\bB}

\def\G{\bG}
\def\T{\bT}

\newcommand{\PP}{\mathbf{P}}

\newcommand{\RR}{\mathbf{R}}
\newcommand{\ZZ}{\mathbf{Z}}

\DeclareMathOperator\val{val}

\DeclareMathOperator\Prin{Prin}

\DeclareMathOperator\Jac{Jac}
\DeclareMathOperator\trop{trop}

\renewcommand\red{\operatorname{red}}

\newcommand\an{{\operatorname{an}}}
\renewcommand\div{{\operatorname{div}}}
\renewcommand\sm{{\operatorname{sm}}}

\newcommand{\p}{{}'}

\newcommand\scdot{\,\cdot\,}

\renewcommand{\deg}{{\operatorname{deg}}}

\title[Lifting meromorphic functions on tropical curves]{The skeleton of the Jacobian, the Jacobian of the skeleton, and lifting meromorphic functions from tropical to algebraic curves}

\author{Matthew Baker} 
\email{mbaker@math.gatech.edu}
\address{School of Mathematics, Georgia Institute of Technology, Atlanta GA 30332-0160, USA}
\author{Joseph Rabinoff}
\email{jrabinoff@math.gatech.edu}
\address{School of Mathematics, Georgia Institute of Technology, Atlanta GA 30332-0160, USA}

\thanks{The authors thank Grisha Mikhalkin for suggesting the ``at most $g$ extra zeros and poles'' assertion in Theorem~\ref{thm:prin.surj}. M.B. was
partially supported by NSF grant DMS-1201473.}

\begin{document}

\begin{abstract}
Let $K$ be an algebraically closed field which is complete with respect to a nontrivial, non-Archimedean valuation and let $\Lambda$ be its value group.  
Given a smooth, proper, connected
$K$-curve $X$ and a skeleton $\Gamma$ of the Berkovich analytification $X^{\an}$, there are two natural real tori which one can consider: the tropical Jacobian 
$\Jac(\Gamma)$ and the skeleton of the Berkovich analytification $\Jac(X)^\an$.  We show that the skeleton of the Jacobian is canonically isomorphic to the
Jacobian of the skeleton as principally polarized tropical abelian varieties.  In addition, we show that the tropicalization of a classical Abel-Jacobi map is a tropical Abel-Jacobi map.  As a consequence of these results,
we deduce that $\Lambda$-rational principal divisors on $\Gamma$, in the
sense of tropical geometry, are exactly the retractions of principal divisors on 
$X$.  We actually prove a more precise result which says that, although zeros and poles of divisors can cancel under the retraction map, in order to lift a $\Lambda$-rational
principal divisor on $\Gamma$ to a principal divisor on $X$ it is never necessary to add more than $g$ extra zeros and $g$ extra poles.
Our results imply that a continuous function $F:\Gamma\to\R$ is the restriction to $\Gamma$ of $-\log|f|$ for some nonzero meromorphic function $f$ on $X$
if and only if $F$ is a $\Lambda$-rational tropical meromorphic function, and we use this fact to prove that there is a
rational map $f : X \dashrightarrow {\mathbf P}^{3}$ whose tropicalization, when restricted to $\Gamma$, is an isometry onto its image.
\end{abstract}

\maketitle

Throughout this paper, $K$ denotes a field which is complete with respect
to a nontrivial, non-Archimedean valuation $\val: K\to\R\cup\{\infty\}$.
Let $R$ be its valuation ring, let $k$ be its residue field, and let 
$\Lambda = \val(K^\times)$ be its value group.  Note that $\Lambda$ is
either discrete or dense in $\R$.  Let $|\scdot| = \exp(-\val(\scdot))$ be
an associated absolute value.

\section{Introduction}
\label{sec:intro}

In this section, we make the additional assumption that $K$ is algebraically closed.
Let $X$ be a smooth, projective, connected $K$-curve and let
$\fX$ be a semistable $R$-model of $X$.  Let $X^\an$ be the
analytification of $X$ in 
the sense of Berkovich~\cite{berkovich:analytic_geometry}.  The incidence graph $\Gamma$ of the
special fiber of $\fX$ naturally has the structure of a finite metric
graph with edge lengths in $\Lambda$; we let $\Gamma(\Lambda)$ denote the
set of all points of $\Gamma$ whose distance from any vertex belongs to $\Lambda$.  
There is a natural inclusion $\Gamma\inject X^\an$
and a deformation retraction $\tau:X^\an\surject\Gamma$.  The metric graph $\Gamma$
is called the \emph{skeleton} of $X$ (or $X^\an$) associated to the model $\fX$.
The retraction $\tau$ takes
$X(K)$ surjectively onto $\Gamma(\Lambda)$, hence defines (by
linearity) a surjective homomorphism 
\[ \tau_* ~:~ \Div(X) \To \Div_\Lambda(\Gamma) \]
of divisor groups, where $\Div_\Lambda(\Gamma)$ is the group of all
$\Z$-linear combinations of points of $\Gamma(\Lambda)$.  Since $\tau_*$ preserves
the degree of a divisor, it also induces a map on degree-zero divisors
$\tau_*:\Div^0(X)\to\Div^0_\Lambda(\Gamma)$.

There is a well-developed theory of divisors and linear equivalence on
graphs and on metric graphs recalled briefly here and in 
\S\ref{section:Tropical Jacobians}; see  \cite{baker:specialization}
and the references therein for details. 
A \emph{tropical meromorphic function} on $\Gamma$ is a continuous function
$F:\Gamma\to\R$ which is piecewise affine with integer slopes.
We say that $F$ is \emph{$\Lambda$-rational}
provided that $F(\Gamma(\Lambda))\subset\Lambda$ and
all points at which $F$ is not differentiable are contained in
$\Gamma(\Lambda)$.   For such an $F$ we define the
\emph{divisor} of $F$ to be 
$\div(F) = \sum_{x\in\Gamma(\Lambda)} n_x (x)$, where $n_x$ is the sum
of the outgoing slopes of $F$ at $x$.  The group of ($\Lambda$-rational)
\emph{principal divisors} on $\Gamma$ is the subgroup
$\Prin_\Lambda(\Gamma)\subset\Div^0_\Lambda(\Gamma)$ of divisors of
$\Lambda$-rational tropical meromorphic functions.

Let $f\in K(X)$ be a nonzero rational function and let $F$ be the
restriction of $-\log|f|$ to $\Gamma\subset X^\an$.  It is a consequence of
non-trivial results in non-Archimedean analysis (see
Remark~\ref{rem:poincare.lelong}) that $F$ is a $\Lambda$-rational
tropical meromorphic function on $\Gamma$ and that 
$\tau_*(\div(f)) = \div(F)$.  In particular, $\tau_*$ takes principal
divisors on $X$ to $\Lambda$-rational principal divisors on $\Gamma$.
One of the main goals of this paper is to prove the following theorem:

\begin{thm} \label{thm:prin.surj}
  Let $K$ be an algebraically closed complete non-Archimedean field, let
  $X$ be a $K$-curve of genus $g$,  
  and let $\Gamma$ be a skeleton of $X$, as above. Then the map on principal divisors
  \[ \tau_* ~:~ \Prin(X) \To \Prin_\Lambda(\Gamma) \]
  is surjective.  More precisely, for $\bar D\in\Prin_\Lambda(\Gamma)$, writing 
  $\bar D = \bar D_1 - \bar D_2$ for effective divisors 
  $\bar D_1,\bar D_2\in\Div^n_\Lambda(\Gamma)$ of degree $n$, there exist
  effective divisors $D_1,D_2\in\Div^{n+g}(X)$ of degree $n+g$ 
  such that $D = D_1-D_2\in\Prin(X)$ and $\tau_*(D) = \bar D$.
\end{thm}

\smallskip
The second assertion states that, although zeros and poles of divisors can
``cancel out'' under the retraction map $\tau_*$, in order to lift a $\Lambda$-rational
principal divisor on $\Gamma$ to a principal divisor on $X$ it is never necessary to add more than $g$ zeros and
poles.  Theorem~\ref{thm:prin.surj} is interesting for several reasons.  From the viewpoint of
tropical geometry, it says that the $\Lambda$-rational principal divisors on $\Gamma$ (which one can think of as an ``abstract complete tropical curve'') are exactly the 
retractions of the principal divisors on $X$ for \emph{any} curve $X$ containing $\Gamma$ as a skeleton.  
(Note that if we replace \emph{any} with \emph{some} then this weaker statement is a consequence of 
\cite[Theorem 9.9]{abbr:hmmc}.)
This gives a precise and very close connection between the notions of principal divisors in tropical and algebraic geometry.

\smallskip

Theorem~\ref{thm:prin.surj} is also interesting from the point of view of
non-Archimedean analysis.  The following corollary is immediate:

\begin{cor} \label{cor:merom.surj}
  In the situation of Theorem~\ref{thm:prin.surj}, 
  let $F:\Gamma\to\R$ be a continuous function.  Then there exists a
  nonzero rational function $f\in K(X)$ such that 
  $F = -\log|f|\,\big|_\Gamma$ if and only if $F$ is a $\Lambda$-rational tropical meromorphic function.
\end{cor}

\smallskip
Among other things, Corollary~\ref{cor:merom.surj} can be used to simplify and sharpen the proofs of some of the ``faithful tropicalization'' results in
\cite{bpr:trop_curves}.  For example, we show in \S\ref{section:faithful_trop} that given a $K$-curve $X$ together with a skeleton $\Gamma$, there is a
rational map $f : X \dashrightarrow {\mathbf P}^{3}$ whose tropicalization, when restricted to $\Gamma$, is an isometry onto its image 
(with respect to lattice lengths in ${\mathbf R}^{3}$).  The main point is that Corollary~\ref{cor:merom.surj} allows us to reduce such questions to purely
combinatorial problems about metric graphs.

\smallskip

The other main goal of this paper is to prove a natural compatibility between classical and tropical Abel-Jacobi maps.
This compatibility forms the backbone of our proof of
Theorem~\ref{thm:prin.surj}, so our two main goals are closely related.
To be more precise, let $X$ be a $K$-curve as above and let $\Gamma$ be a skeleton for $X$.
Choose a $K$-point $P \in X(K)$ and let $p = \tau(P)\in\Gamma(\Lambda)$.  Let $\alpha_P: X\to J$ and $\alpha_p: \Gamma\to\Jac(\Gamma)$ be the classical and tropical Abel-Jacobi maps based at $P$ and $p$, respectively.  (The tropical Abel-Jacobi map is defined in \S\ref{section:Tropical Jacobians} below.)   We prove:
 
\begin{thm} \label{thm:abel-jac_compat}  
There is a canonical isomorphism $\Jac(\Gamma) \cong \Sigma(J^\an)$, where $\Sigma(J^\an)$ is the skeleton of $J^\an$ in the sense of
 Berkovich, such that (identifying $\alpha_p$ with the corresponding
 map from $\Gamma$ to $\Sigma(J^\an)$) the following square commutes: 
  \begin{equation} \label{eq:intro.aj.commutes}\xymatrix @=.25in{
       {X^\an} \ar[r]^{\alpha_P} \ar[d] & 
      {J^\an} \ar[d] \\
      {\Gamma} \ar[r]_(.4){\alpha_p} & {\Sigma(J^\an)}
    }\end{equation}
where the vertical arrows are the natural retraction maps.
Moreover, the isomorphism $\Jac(\Gamma) \cong \Sigma(J^\an)$ is compatible with the natural principal polarizations on both sides (so defines
an isomorphism of {\em principally polarized tropical abelian varieties}).
\end{thm}

\smallskip

See Corollary~\ref{cor:skeleta.compatibility} and and
Proposition~\ref{prop:compat.aj}. 
The fact that $\Jac(\Gamma) \cong \Sigma(J^\an)$ can be summarized as:
\begin{center}
\emph{The skeleton of the Jacobian is the Jacobian of the skeleton}. 
\end{center}

\smallskip

A closely related result has recently been established in \cite[Theorem A]{viviani:trop_vs_compact}.

\section{The analogy with arithmetic geometry; strategy of proof}

In order to explain how our two main goals are related, and to outline the
proofs of Theorems~\ref{thm:prin.surj} and \ref{thm:abel-jac_compat}, it
is useful to make a digression into what was previously known in the
analogous situation where $K$ is discretely valued.  This analogy with
arithmetic geometry was one of the original motivations for this study.

\paragraph[The discretely valued situation] \label{par:discrete.story}
Suppose that the 
valuation on $K$ is discrete with value group $\Lambda=\Z$ and that 
the residue field $k$ is algebraically closed.  
As in \S\ref{sec:intro} we let $X$ be a smooth, proper, geometrically connected $K$-curve
with semistable $R$-model $\fX$.
For simplicity we also assume that $\fX$ is regular
and that $\fX_k$ is (split) semistable with smooth irreducible components.
In this case $\Gamma(\Lambda) = \Gamma(\Z)$ is the set of
vertices of $\Gamma$, i.e.\ all edges have length one, so we may regard
$\Gamma$ as an ordinary (non-metric) graph $G$; then the theory of $\Lambda$-rational divisors
and linear equivalence on $\Gamma$ reduces to Baker and Norine's notions of divisors and linear equivalence on $G$
from~\cite{baker_norine:rr_graphs}.  

The abelian group 
\[ \Jac(G) = \Jac_\Z(\Gamma) \coloneq \Div^0_\Z(\Gamma)/\Prin_\Z(\Gamma) \]
is called the \emph{Jacobian} of $G$ (or the \emph{regularized Jacobian} of $\Gamma$). 
It has the following alternative description.  Let $M = H_1(\Gamma,\Z)$ be the first homology
group of (the topological realization of) $\Gamma$ and let
$N = \Hom(M,\Z)$ be its dual.  We regard $M$ as the
group of cycles in the free abelian group $C_1(G,\Z)$ generated by
the edges of $G$.  Define a pairing $\angles{\scdot,\scdot}$ on
$C_1(G,\Z)$ by
\begin{equation} \label{eq:edge.length.pairing.1}
  \angles{e,e'} =
  \begin{cases}
    1 &\quad\text{ if } e = e' \\
    0 &\quad\text{ otherwise.}
  \end{cases}
\end{equation}
This restricts to a symmetric, positive-definite pairing 
$\angles{\scdot,\scdot}:M\times M\to\Z$, hence defines an injective
homomorphism $\mu: M\to N$.

\begin{prop} \label{prop:discrete.div.unif.jac}
  There is a canonical short exact sequence
  \[ 0 \To M \overset\mu\To N \To \Jac(G) \To 0. \]
\end{prop}

\smallskip
Proposition~\ref{prop:discrete.div.unif.jac} should be viewed as a
compatibility between the divisor-theoretic description of
$\Jac_\Z(\Gamma)$, and a description in terms of 
``the dual of a space of one-forms modulo a period lattice.'' 
See \S\ref{section:Tropical Jacobians} for details.

Let $J$ be the Jacobian of $X$, let $\fJ$ be its
N\'eron model over $R$, let $\fJ^0$ be the (fiberwise) connected component
of the identity in $\fJ$, and let $\Phi = \fJ_k/\fJ^0_k = J(K)/\fJ^0(R)$ be the
component group (the second equality holds
using~\cite[Proposition~10.1.40(b)]{liu:algebraic_geometry}).  Since $X$
has split semistable reduction, it is known that $\fJ^0_k$ is a split
semi-abelian variety.

By a general theorem of Grothendieck 
\cite[Expos\'e~IX, 11.5]{SGA7}, the component group $\Phi$ fits into an exact sequence of the form
  \[ 0 \To M' \xrightarrow{\mu_{\rm mon}} N' \To \Phi \To 0, \]
where $M'$ is the character lattice of the toric part of $\fJ^0_k$, $N'$ is its dual, and 
$\mu_{\rm mon}$ is the natural map derived from the {\em monodromy pairing}
$\angles{\scdot,\scdot}_{\rm mon}:M'\times M'\to\Z$.

A proof of the following result can be found in \cite{bosch_lorenzini:monodromy}; see
also~\cite{SGA7,illusie:monodromy.pairing}.

\begin{thm} \label{thm:analytic.monodromy.pairing}
  The character lattice $M'$ of the toric part of $\fJ^0_k$ is canonically isomorphic
  to $M$, and the pairing $\angles{\scdot,\scdot}$ defined above coincides
  with Grothendieck's monodromy pairing under this isomorphism.
\end{thm}

\smallskip

The previously published proofs of this result use arithmetic intersection
theory and rest on a Picard-Lefschetz argument.  We are able to
derive (the second part of) Theorem~\ref{thm:analytic.monodromy.pairing}
from our results, thus giving a very different proof; see
Remark~\ref{rem:second.proof}.

As an immediate consequence of Proposition~\ref{prop:discrete.div.unif.jac} and Theorem~\ref{thm:analytic.monodromy.pairing}, one has the following 
discrete analogue of the second part of Theorem~\ref{thm:abel-jac_compat}:

\begin{cor} \label{cor:comp.gp.monodromy.pairing}
  There is a canonical isomorphism 
  $\Jac(G)\cong\Phi$.
\end{cor}

\smallskip

Corollary~\ref{cor:comp.gp.monodromy.pairing} is a reformulation of Raynaud's famous description of the component group $\Phi$ in terms of
the intersection matrix of the special fiber of $\fX$ \cite[Proposition 8.1.2]{raynaud:component_groups}, since the latter coincides with the Laplacian matrix of $G$.
Raynaud's result plays an important technical role in many key papers on the arithmetic of modular curves, including for example \cite{mazur:eisenstein.ideal, 
ribet:modular.repns}.

We next turn to a discrete analogue of Theorem~\ref{thm:prin.surj}.  For this, we define the
{\em specialization map} $\rho: \Div(X)\to\Div_\Z(\Gamma)$ from divisors on $X$ to divisors on $G$ as follows.  
For $x\in\Gamma(\Z)$ let $C_x$ be the associated irreducible component of
$\fX_k$, considered as a divisor on $\fX$.  If $D\in\Div(X)$ is a divisor
then its Zariski closure $\fD$ in $\fX$ is again a Weil divisor; we set
\[ \rho(D) = \sum_{x\in\Gamma(\Z)} (C_x\cdot\fD)\,(x), \]
where $C_x\cdot\fD$ is the intersection number
$\deg(\sO_\fX(\fD)|_{C_x})$.  (Note that this map $\rho$ does \emph{not}
coincide with the extension by linearity of the retraction map 
$\tau: X(\bar K)\to \Gamma$, as the image of the latter is not contained in
$\Gamma(\Z)$.  However, $\tau_*(D)$ and $\rho(D)$ are always linearly
equivalent divisors: see \cite[Remark 2.10]{baker:specialization}.)  
By \cite[Lemma 2.1]{baker:specialization},
$\rho$ takes principal divisors to principal divisors, i.e.,\ restricts to a map
$\rho_0:\Prin(X)\to\Prin_\Z(\Gamma)$.

\begin{thm} \label{thm:discrete.retract.ses} 
  The following diagram
  \begin{equation} \label{eq:discrete.retract.ses} 
    \xymatrix @R=.25in{
      0 \ar[r] & 
      {\Prin(X)} \ar[r] \ar[d]^{\rho_0} &
      {\Div^0(X)} \ar[r] \ar[d]^\rho &
      {J(K)} \ar[r] \ar[d] & 0 \\
      0 \ar[r] &
      {\Prin_\Z(\Gamma)} \ar[r] &
      {\Div^0_\Z(\Gamma)} \ar[r] & 
      {\Jac_\Z(\Gamma)=\Phi} \ar[r] & 0
    }\end{equation}
  is commutative, where the right vertical arrow is the canonical quotient
  map $J(K)\to J(K)/\fJ^0(R)=\Phi$.  
\end{thm}

\smallskip
Theorem~\ref{thm:discrete.retract.ses} is a consequence of a theorem of
Raynaud~\cite[Theorem~9.6.1]{blr:neron_models}, as was observed by the first
author~\cite[Appendix~A]{baker:specialization}.  

The specialization map from $X(K)$ to the smooth locus of $\fX_k(k)$ is
surjective (see for
instance \cite[Proposition~10.1.40(b)]{liu:algebraic_geometry}); hence
$\rho:\Div^0(X)\to\Div^0_\Z(\Gamma)$ is surjective.  Applying the
snake lemma to the morphism of short exact
sequences~\eqref{eq:discrete.retract.ses}, we obtain an exact sequence
\[ 0 \To \Prin^{(0)}(X) \To \Div^{(0)}(X) \To \fJ^0(R) \To \coker(\rho_0)
\To 0, \]
where $\Div^{(0)}(X) = \ker(\rho)$ and 
$\Prin^{(0)}(X) = \ker(\rho_0) = \Div^{(0)}(X)\cap \Prin(X)$.
A celebrated theorem of
Raynaud~\cite{raynaud:component_groups}, \cite[Theorem~9.5.4]{blr:neron_models} 
asserts that $\fJ^0$ is the relative Jacobian of $\fX$.  That is to say,
$\fJ^0$ represents the functor of line bundles of total degree zero, which
implies that $\fJ^0(R)$ is the group of those line bundles on $J$
represented by divisors $D$ such that $\rho(D) = 0$.  Therefore we have 
$\fJ^0(R) = \Div^{(0)}(X)/\Prin^{(0)}(X)$, so we obtain the
discretely-valued version of Theorem~\ref{thm:prin.surj}:

\begin{thm}[{\cite[Appendix~A]{baker:specialization}}]
  \label{thm:discrete.prin.surj}
  With the above hypotheses, the specialization map
  \[ \rho_0 ~: ~ \Prin(X)\To\Prin_\Z(\Gamma) \]
  is surjective. 
\end{thm}

\paragraph[Strategy of the proofs]
Much of this paper is devoted to developing analogues 
of the tools and the language used above for more general 
(not necessarily discretely valued) non-Archimedean fields $K$.
Our ``tropical'' point of view on these results forms an attractive framework which should be useful
in its own right.
The primary difficulty is that N\'eron models are not available 
in the non-Noetherian situation.  Instead we use the non-Archimedean
uniformization theory of Jacobians, which was worked out by
Bosch, L\"utkebohmert, and Raynaud in this generality.

Let $K$ be any complete non-Archimedean field and let
$X$ be a smooth, projective, geometrically connected $K$-curve endowed
with a split semistable $R$-model $\fX$.  As above we let 
$\Gamma\subset X^\an$ be the associated skeleton, with retraction map
$\tau:X^\an\to\Gamma$.  The $R$-model $\fX$ defines a weighted graph model $G$ for $\Gamma$ in the sense of (\ref{par:graph.basics}) below;
for ease of terminology we will refer to vertices and edges of $\Gamma$ instead of $G$.

Recall that $\Div^0_\Lambda(\Gamma)$ is the group of degree-zero divisors
on $\Gamma$ supported on $\Gamma(\Lambda)$ and that
$\Prin_\Lambda(\Gamma)$ is the 
group of divisors of $\Lambda$-rational tropical meromorphic functions on
$\Gamma$.  The \emph{set of $\Lambda$-points of the Jacobian of $\Gamma$}
is defined to be the abelian group
\[ \Jac_\Lambda(\Gamma) = \Div^0_\Lambda(\Gamma) / \Prin_\Lambda(\Gamma). \]
This group has the following alternative description.  Let
$M = H_1(\Gamma,\Z)$ and let $N = \Hom(M,\Z)$.
The \emph{edge length pairing} $\angles{\scdot,\scdot}$ on the free
abelian group $C_1(\Gamma,\Z)$ generated by the edges of $\Gamma$ is
defined on edges by  
\begin{equation} \label{eq:edge.length.pairing}
  \angles{e,e'} =
  \begin{cases}
    \ell(e) &\quad\text{ if } e = e' \\
    0 &\quad\text{ otherwise,}
  \end{cases}
\end{equation}
where $\ell(e)$ is the length of $e$.  Regarding 
$M = H_1(\Gamma,\Z)$ as a subgroup of $C_1(\Gamma,\Z)$, 
we restrict $\angles{\scdot,\scdot}$ to a symmetric, positive-definite
pairing on $M$.  This is exactly analogous
to~\eqref{eq:edge.length.pairing.1}.  Since $\angles{\scdot,\scdot}$ takes
values in $\Lambda$, we obtain a map 
$\mu: M\to N_\Lambda = \Hom(M,\Lambda)\subset N_\R$.
We have the following metric graph generalization of
Proposition~\ref{prop:discrete.div.unif.jac}: 

\begin{prop} \label{prop:div.unif.jac}
  There is a canonical short exact sequence
  \[ 0 \To M \overset\mu\To N_\Lambda \To \Jac_\Lambda(\Gamma) \To 0. \]
\end{prop}

\smallskip
See Proposition~\ref{prop:div.unif.jac2}.
We define the \emph{Jacobian} of $\Gamma$ to be the real torus
$\Jac(\Gamma) \coloneq \Jac_\R(\Gamma) = N_\R / \mu(M)$.
Note that $\Jac_\Lambda(\Gamma)\subset\Jac(\Gamma)$.

Let $J$ be the Jacobian of $X$ and let $E$ be the
universal cover of $J^\an$.  This topological space has the structure of
a $K$-analytic group making the structure map $E\to J^\an$ into a
$K$-analytic group homomorphism.  Its kernel 
$H_1(J^\an,\Z)\subset E(K)$ is canonically isomorphic to
$M = H_1(\Gamma,\Z) = H_1(X^\an,\Z)$.   There exists a  
continuous map $\trop:E\to  N_\R = \Hom(M,\R)$, which restricts to an
isomorphism $\trop: M\isom\trop(M)$ onto a full-rank lattice
$\trop(M)\subset N_\R$.  Letting $\Sigma$ be the real torus 
$N_\R/\trop(M)$, we obtain a commutative diagram
\[\xymatrix @R=.2in{
  0 \ar[r] & M \ar[r] \ar[d]^\trop & E \ar[r] \ar[d]^\trop & 
  {J^\an} \ar[r] \ar[d]^{\bar\tau} & 0 \\
  0 \ar[r] & {\trop(M)} \ar[r] & {N_\R} \ar[r] & 
  \Sigma \ar[r] & 0 
}\]
The real torus $\Sigma$ is called the \emph{skeleton} of the abelian
variety $J$.  Berkovich has
shown~\cite[\S6.5]{berkovich:analytic_geometry} that there is an inclusion 
$\Sigma\inject J^\an$, and that $\bar\tau: J^\an\to\Sigma$ is a
deformation retraction.  

The map $\trop: M\inject E^\an\to N_\R$ defines a bilinear pairing $\angles{\scdot,\scdot}_\an:M\times M\to\R$
which makes the triple $(\Sigma,M,\angles{\scdot,\scdot}_\an)$ into a
principally polarized tropical abelian variety in the sense of~\parref{par:pptav}.
One of the main results of this paper is the following analogue of Theorem~\ref{thm:analytic.monodromy.pairing}.

\begin{thm} \label{thm:skeleta.compatibility}
  The pairings $\angles{\scdot,\scdot}_\an$ and $\angles{\scdot,\scdot}$
  coincide.  As a consequence, there is a canonical isomorphism of 
  principally polarized tropical abelian varieties $\Sigma\cong\Jac(\Gamma)$.
\end{thm}

\smallskip

In the discretely-valued case, Theorem~\ref{thm:skeleta.compatibility} is exactly
Theorem~\ref{thm:analytic.monodromy.pairing}.  In fact, it is easy to
derive Theorem~\ref{thm:analytic.monodromy.pairing} from
Theorem~\ref{thm:skeleta.compatibility}, thus giving a very different
proof; see Remark~\ref{rem:second.proof}.  
In the special case where $X$ is a Mumford curve, 
Theorem~\ref{thm:skeleta.compatibility} has also been proved by
van~der~Put~\cite[Theorem~6.4]{vdPut:discrete.mumford.theta},
using the theory of $p$-adic theta functions.
Our proof in the general case is based on Theorem~\ref{thm:abel-jac_compat} (the compatibility between
the algebraic Abel-Jacobi map and its tropical analogue).

\begin{rem} \label{rem:trop.of.aj}
  Theorem~\ref{thm:skeleta.compatibility} has the following 
  interpretation.  Let $g\geq 2$.  There is a map
  $\trop: M_g^\an\to\Djunion_{g'\leq g} M_{g'}^{\trop}$ from the
  analytification of the moduli space of genus $g$ curves to the moduli
  space of metric graphs of genus $g'\leq g$, which takes a curve $X$ to
  its minimal skeleton.  There is also a 
  map $\trop: A_g^\an\to\Djunion_{g'\leq g} A_{g'}^{\trop}$ from the moduli
  space of principally polarized abelian varieties of genus $g$ to the
  moduli space of principally polarized tropical abelian varieties of dimension $g'\leq g$,
  taking an abelian variety to its skeleton in the sense of Berkovich.
  There are Torelli maps $M_g\to A_g$ and $M_{g'}^{\trop}\to A_{g'}^{\trop}$
  which take a curve (resp.\ graph) to its Jacobian. 
  Theorem~\ref{thm:skeleta.compatibility} is equivalent to the statement that the following square
  commutes: 
  \[\xymatrix @R=.25in{
    {M_g^\an} \ar[r]^(.45){\trop} \ar[d] & {M_g^{\trop}} \ar[d] \\
    {A_g^\an} \ar[r]^(.45){\trop} & {A_g^{\trop}}
  }\]
  See also \cite[Theorem A]{viviani:trop_vs_compact}, and compare with Remark~\ref{rem:period.matrix}.
\end{rem}

\smallskip
The retraction map $\tau: X(K)\to\Gamma(\Lambda)$
extends by linearity to give a map
$\tau_*:\Div^0_K(X)\to\Div^0_\Lambda(\Gamma)$, where
$\Div^0_K(X)$ is the group of degree-zero divisors on $X$ supported on
$X(K)$.  Let $\Prin_K(X)$ be the kernel of the natural homomorphism
$\Div^0_K(X)\to J(K)$.  One checks
after passing to the completion of the algebraic closure of $K$ that 
$\tau_*(\Prin_K(X))\subset\Prin_\Lambda(\Gamma)$.
We will prove the following analogue of
Theorem~\ref{thm:discrete.retract.ses}, which says that
the retraction map $\bar\tau: J^\an\to\Sigma = \Jac(\Gamma)$ is compatible
with the descriptions of $J(K)$ and $\Jac_\Lambda(\Gamma)$  in terms of
divisors modulo principal divisors:

\begin{thm} \label{thm:retract.ses}
  We have $\bar\tau(J(K))\subset\Jac_\Lambda(\Gamma)$,
  and the following diagram is commutative:
  \[\xymatrix @R=.25in{
    0 \ar[r] & {\Prin_K(X)} \ar[r] \ar[d]^{\tau_*}
    & {\Div^0_K(X)}  \ar[r] \ar[d]^{\tau_*} 
    & {J(K)} \ar[d]^{\bar\tau} &  \\
    0 \ar[r] & {\Prin_\Lambda(\Gamma)} \ar[r] 
    & {\Div^0_\Lambda(\Gamma)} \ar[r]
    & {\Jac_\Lambda(\Gamma)} \ar[r] & 0
  }\]
\end{thm}

\smallskip
Note that the arrow $\Div_K^0(X)\to J(K)$ is surjective when $K$ is
algebraically closed.
The first part of Theorem~\ref{thm:prin.surj} follows from
Theorem~\ref{thm:retract.ses} in 
the same way that Theorem~\ref{thm:discrete.prin.surj} follows from
Theorem~\ref{thm:discrete.retract.ses} (applying the snake lemma to the above diagram), with a result of Bosch and
L\"utkebohmert playing the role of Raynaud's theorem that the connected
component of the N\'eron model of $J$ represents the relative Jacobian.
The second assertion in Theorem~\ref{thm:prin.surj} requires a separate
argument.

\begin{rem}
  In fact we will prove a preliminary version of
  Theorem~\ref{thm:retract.ses} \emph{before} proving 
  Theorem~\ref{thm:skeleta.compatibility}: see
  Proposition~\ref{prop:define.taubar}.  This is necessary because 
  Theorem~\ref{thm:skeleta.compatibility} is proven using a compatibility
  of Abel-Jacobi maps, which already depends on the divisor-theoretic
  description of the retraction map $\bar\tau$.  
\end{rem}

\section{The Jacobian of a metric $\Lambda$-graph}
\label{section:Tropical Jacobians}

We briefly review the theory of metric graphs and their Jacobians,
following \cite{MZ} (see also \cite{BF}). 

\paragraph \label{par:graph.basics}
Let $\Lambda$ be a nonzero additive subgroup of $\RR$ and let $\Gamma$ be a $\Lambda$-metric graph in the sense of \cite[\S{1.2}]{abbr:hmmc}.
A $\Lambda$-rational {\em model} $\Gamma$ is a weighted graph $G$ with edge lengths in $\Lambda$ whose geometric realization is $\Gamma$; see \cite[\S{1.2}]{abbr:hmmc} or \cite[\S{3}]{BF}.
Recall that $\Gamma(\Lambda)$ denotes the subset of $\Lambda$-rational points of $\Gamma$, which formally means the set of points of $\Gamma$ whose
distance to any vertex of $G$ belongs to $\Lambda$, for any $\Lambda$-rational model $G$ for $\Gamma$.
We define $\Div_\Lambda(\Gamma)$ to be the free abelian group on $\Gamma(\Lambda)$ and
$\Div^0_\Lambda(\Gamma)$ to be its degree-zero subgroup.
The group $\Prin_\Lambda(\Gamma)$ consisting of divisors of $\Lambda$-rational tropical meromorphic functions on
$\Gamma$ is a subgroup of $\Div^0_\Lambda(\Gamma)$, and we define $\Pic^0_\Lambda(\Gamma)$ to be the quotient group:
\[
\Pic^0_\Lambda(\Gamma) = \Div^0_\Lambda(\Gamma) / \Prin_\Lambda(\Gamma).
\]
We also set $\Pic^0(\Gamma) = \Pic^0_\R(\Gamma)$.

\smallskip

Following \cite[Definition 4.1]{MZ} and \cite[\S{2.1},\S{3}]{BF}, we define the group $\Omega_{\Gamma(\Lambda)}(\Gamma)$ of {\em $\Lambda$-harmonic $1$-forms} on $\Gamma$ to be the direct limit of $\Omega(G)$ over all $\Lambda$-rational weighted graph models $G$ for $\Gamma$.
(This coincides with the subgroup of $\Lambda$-rational $1$-forms in the sense of Mikhalkin-Zharkov by \cite[Footnote~5]{BF}.)
We can view an element $\gamma \in H_1(\Gamma,\ZZ)$ as an element of 
\[
\Omega_{\Gamma(\Lambda)}(\Gamma)^* = \Hom\big(\Omega_{\Gamma(\Lambda)}(\Gamma),\Lambda\big)
\]
via the integration pairing $\mu: \gamma \mapsto \int_\gamma$.
We define 
\[
\Jac_\Lambda(\Gamma) = \Omega_{\Gamma(\Lambda)}(\Gamma)^* /
\mu(H_1(\Gamma,\ZZ))
\sptxt{and} \Jac(\Gamma) = \Jac_\R(\Gamma).
\]
If we set $M=H_1(G,\ZZ)$ then it is formal to check that if we identify $\Omega_{\Gamma(\Lambda)}(\Gamma)$ with $H_1(\Gamma,\Lambda)$
in the obvious way and set $N_{\Lambda}=\Hom(M,\Lambda)$ then the above definition becomes
\[
\Jac_\Lambda(\Gamma) = N_{\Lambda} / \mu(H_1(\Gamma,\ZZ)),
\]
where $\mu$ now denotes the map from $H_1(\Gamma,\ZZ)$ to $N_{\Lambda}$
given by the edge length pairing~\eqref{eq:edge.length.pairing}.

The following result, a reformulation of Proposition~\ref{prop:div.unif.jac}, was first proved for $\Lambda = \RR$ in \cite[Theorem 6.2]{MZ}.  
An alternate proof, which easily generalizes to the present situation in which $\Lambda$ is an arbitrary subgroup of $\RR$, is given in \cite[Theorem 3.4]{BF}.

\begin{prop} \label{prop:div.unif.jac2}
There is a canonical isomorphism $\Pic^0_\Lambda(\Gamma) \cong \Jac_\Lambda(\Gamma)$.
\end{prop}

\smallskip
In particular, when $\Lambda = \RR$ there is a canonical isomorphism between the degree-zero Picard group $\Pic^0(\Gamma)$ and the $g$-dimensional
real torus $\Jac(\Gamma)$, where $g = \dim_{\RR} H^1(\Gamma,\RR)$.

\paragraph \label{par:tropical.aj}
If we fix a base point $q \in \Gamma$, we obtain a corresponding (tropical) {\em Abel-Jacobi map}
\[
\alpha = \alpha_q : \Gamma \to \Jac(\Gamma)
\]
taking $\Gamma(\Lambda)$ to $\Jac_{\Lambda}(\Gamma)$ if $q \in \Gamma(\Lambda)$.
The map $\alpha$ is defined by 
\[
\alpha({p}) = \int_q^p \in \Omega_{\Gamma(\R)}(\Gamma)^*,
\]
where the integral is computed along any path in $\Gamma$ from $q$ to $p$.
Different choices of paths give rise to elements of $\Omega_{\Gamma(\R)}(\Gamma)^*$ which differ by an element of $\mu(H_1(\Gamma,\ZZ))$.
From this description, it is straightforward to check (just as in the classical case of Riemann surfaces) that we have:

\begin{lem} \label{lem:graph.aj.kpi1}
  The Abel-Jacobi map $\alpha : \Gamma \to \Jac(\Gamma)$ induces an isomorphism 
  \[ \alpha_* : H_1(\Gamma,\ZZ) \isom H_1(\Jac(\Gamma),\ZZ) \]
  on singular homology groups.
\end{lem}

\smallskip
We remark that under the canonical isomorphism afforded by Proposition~\ref{prop:div.unif.jac2}, 
the map $\alpha$ corresponds to the map $\Gamma \to \Pic^0(\Gamma)$ given by $p \mapsto [({p)}-(q)]$.

\paragraph
As in classical algebraic geometry, there is a canonical bijection between isomorphism classes of $\Lambda$-line bundles on $\Gamma$ and
elements of $\Pic_{\Lambda}(\Gamma)$.  Following 
\cite[Definition~4.4]{MZ}, a {\em line bundle} on $\Gamma$ is by
definition a fiber bundle  
$\pi : L \to \Gamma$ whose fibers are isomorphic to the tropical affine line $\T$ and whose transition maps are {\em harmonic} $\Lambda$-rational tropical meromorphic functions on $\Gamma$.
Isomorphism classes of $\Lambda$-line bundles on $\Gamma$ are classified by the \v Cech cohomology group $H^1(\Gamma, {\mathcal O^*_{\Lambda}})$,
where ${\mathcal O^*_{\Lambda}}$ is the sheaf of {\em harmonic} $\Lambda$-rational tropical meromorphic functions on $\Gamma$.

Given an open subset $U \subset \Gamma$ which locally trivializes $L$, one
says that a section $s : U \to \pi^{-1}(U)$ is {\em regular} (resp.\ {\em rational})
if for any open subset $V \subset U$ and any trivialization $\Phi : \pi^{-1}(V) \cong V \times \T$, $s$ becomes a regular (resp. rational) function on $V$.
A global rational section of $L$ gives rise to a well-defined divisor on $\Gamma$, independent of the choice of trivialization.
By \cite[Proposition 4.6]{MZ}, every divisor defines a line bundle together with a rational section which is well-defined up to an additive (i.e., tropically multiplicative) constant, every line bundle has a rational section, and the divisors of any two rational sections are linearly equivalent.

If $\{U_i \}$ is an open cover of $\Gamma$ trivializing $L$ then the map sending $[D]$ to the $1$-cocycle whose value on $U_i$ is (some choice of) a local section $s_i$ associated to $D$ induces a natural map from $\Pic_{\Lambda}(\Gamma)$ to $H^1(\Gamma, {\mathcal O^*_{\Lambda}})$.
By the proof of \cite[Proposition 4.6]{MZ}, we have:

\begin{prop} \label{prop:pic.from.H1}
The natural map
\[
\Pic_{\Lambda}(\Gamma) \to H^1(\Gamma, {\mathcal O^*_{\Lambda}})
\]
is an isomorphism.  
\end{prop}

\paragraph \label{par:pptav}
A {\em principally polarized tropical abelian variety} $A$ of dimension $g$ is a triple $(A,L,\angles{\scdot,\scdot})$, where $A$ is a $g$-dimensional real torus,
$L$ is a rank-$g$ lattice, and $\angles{\scdot,\scdot} : L \times L \to \RR$ is a bilinear form whose associated quadratic form is positive definite,
together with an isomorphism $A \cong \Hom(L,\RR)/\mu(L)$ where $\mu : L \to \Hom(L,\RR)$ is the natural map derived from the bilinear form.
(This is more or less equivalent to the definition given in \cite[\S{5}]{MZ}, and slightly less general then the one used in \cite[Definition 3.1.1]{viviani:trop_vs_compact} and
the references therein.)
An isomorphism $(A,L,\angles{\scdot,\scdot}) \to (A',L',\angles{\scdot,\scdot}')$ of principally polarized tropical abelian varieties is an isomorphism $L \to L'$ which transforms
the bilinear form on $L$ to the one on $L'$, together with the isomorphism $A \to A'$ induced by the resulting identification of $\Hom(L,\RR)/\mu(L)$ with
$\Hom(L',\RR)/\mu'(L')$.

\smallskip

There is a canonical principal polarization on the Jacobian $\Jac(\Gamma)$ of a metric graph $\Gamma$ defined by taking
$L = H_1(\Gamma,\ZZ)$ and $\angles{\scdot,\scdot}$ to be the edge length pairing (\ref{eq:edge.length.pairing}).

\section{Non-Archimedean uniformization of abelian varieties}
\label{sec:general.uniformization}

Let $K$ be any complete non-Archimedean field.  In this section we recall
some results of Bosch, Bosch-L\"utkebohmert, and Raynaud, in Berkovich's language;
see~\cite{bosch:formal_groups,bosch_lutk:uniformizationII,bosch_lutk:degenerating_avs}
and~\cite[\S6.5]{berkovich:analytic_geometry}.  We also prove a
compatibility result involving the Abel-Jacobi map.

\paragraph \label{par:unif.basics}
Let $A$ be an abelian variety over $K$ with split semi-abelian reduction.
This means that there is a unique compact analytic domain 
$A_0\subset A^\an$ which is a formal $K$-analytic subgroup in the sense
of~\cite{bosch:formal_groups}, whose special fiber $\bar A_0$ is an
extension of an abelian variety $\bar B$ by a split torus $\bar T$.
Such an abelian variety admits a non-Archimedean uniformization, in the
following sense.  Let $E^\an$ be the topological universal cover of 
$A^\an$.  Then $E^\an$ has a unique structure of $K$-analytic group (a group
object in the category of $K$-analytic spaces) such that the structure map
$\pi: E^\an\to A^\an$ is a homomorphism of $K$-analytic groups.  
In fact $E^\an$ is the analytification of an algebraic group $E$, although
$\pi$ is not algebraic.  As $\pi$ is a
local isomorphism, its kernel 
$M' \coloneq \ker(\pi)\cong H_1(A^\an,\Z)\cong\pi_1(A^\an)$
is a discrete subgroup of $E^\an(K)$, so we have an exact sequence of $K$-analytic groups
\begin{equation} \label{eq:M.E.A}
  0 \To M' \To E^\an \overset\pi\To A^\an \To 0.
\end{equation}
The algebraic group $E$ is an extension of an abelian variety $B$ with
good reduction by a split torus $T$:
\begin{equation} \label{eq:T.E.B}
  0 \To T \To E \To B \To 0.
\end{equation}
More precisely, $B$ has an abelian scheme $R$-model $\fB$ with special fiber $\bar B$,
and $T$ is the generic fiber of the unique $R$-torus $\fT$ with special
fiber $\bar T$ (so $T$ has the same character lattice as $\bar T$).  Let
$T_0$ denote the affinoid torus inside of 
$T^\an$, i.e.\ the locus of all points in $T^\an$ with a well-defined
specialization in $\bar T$. The short exact sequnce
$0\to\bar T\to\bar A_0\to\bar B\to 0$ lifts to a short exact sequence
\begin{equation} \label{eq:T0.A0.B}
  0 \To T_0 \To A_0 \To B^\an \To 0,
\end{equation}
and~\eqref{eq:T.E.B} is the (algebraization of the) push-out
of~\eqref{eq:T0.A0.B} with respect to the inclusion $T_0\inject T^\an$. 
In particular, $A_0$ is an analytic domain in $E^\an$.
The sequence~\eqref{eq:T.E.B} is unique since $T$ is the maximal torus
in $E$ (the image of any morphism $\G_m\to E$ is contained in
$T$).  The above uniformization data~\eqref{eq:M.E.A}, \eqref{eq:T.E.B}
are often combined into the following diagram, known as the
{\em Raynaud uniformization cross}:
\begin{equation} \label{eq:raynaud.cross}
  \xymatrix @=.20in{
    & {T^\an} \ar[d] & \\
    {M'} \ar[r] & {E^\an} \ar[r] \ar[d] & {A^\an} \\
    & {B^\an} &
  }
\end{equation}
The formation of~\eqref{eq:raynaud.cross} is compatible with extension of
the ground field in the obvious way.

\paragraph \label{par:unif.basics.2}
Let $M$ be the character lattice of $T$ and let $N = \Hom(M,\Z)$, with
$(\scdot,\scdot): N\times M\to\Z$ the evaluation pairing.
For $u\in M$ we let $\chi^u\in K[M]$ denote the corresponding character of
$T$.  We have a tropicalization map
$\trop: T^\an\to N_\R = \Hom(M,\R)$ defined by
$(\trop(\|\scdot\|),\,u) = -\log\|\chi^u\|$, where we regard $T^\an$ as a
space of semi-norms on $K[M]$.  The map $\trop$ is surjective, continuous,
and proper; in fact the affinoid torus $T_0$ is equal to $\trop\inv(0)$.
We extend $\trop$ to all of $E^\an$ by declaring that
$\trop\inv(0) = A_0$, which makes sense because
\[\xymatrix @=.2in {
  {T_0} \ar[r] \ar[d] & {A_0} \ar[d] \\
  {T^\an} \ar[r] & {E^\an}
}\]
is a push-out square.  If $K'$ is any valued field extension of
$K$ with value group $\Lambda'$ then 
$\trop(E(K'))\subset N_{\Lambda'} = \Hom(M,\Lambda')$, and the map
$\trop: E(K')\to N_{\Lambda'}$ is a surjective group homomorphism.
The restriction of $\trop$ to $M'\subset E(K)$ is injective, and
its image $\trop(M')\subset N_\Lambda$ is a full-rank lattice in the real
vector space $N_\R$.  Let $\Sigma = \Sigma(A)$ be the real torus
$N_\R/\trop(M)$.  Since $A^\an$ is the quotient of $E^\an$ by the
action of $M'$, there exists a unique map $\bar\tau: A^\an\to\Sigma$
making the following diagram commute:
\begin{equation} \label{eq:trop.of.unif}
  \xymatrix @R=.2in{
    0 \ar[r] & {M'} \ar[r] \ar[d]^{\trop}_\cong
    & {E^\an} \ar[r] \ar[d]^{\trop}
    & {A^\an} \ar[r] \ar[d]^{\bar\tau} & 0 \\
    0 \ar[r] & {\trop(M')} \ar[r] 
    & {N_\R} \ar[r] & {\Sigma} \ar[r] & 0
  }
\end{equation}
Berkovich~\cite[Theorem~6.5.1]{berkovich:analytic_geometry} has shown that
there exists a section  
$\iota:\Sigma\inject A^\an$ of $\bar\tau$, and that $\iota\circ\bar\tau$
is a deformation retraction.  For this 
reason $\Sigma$ is called the \emph{skeleton} of $A$.
Letting $\Sigma_\Lambda = N_\Lambda/\trop(M)$ and taking $K$-points, we
have a surjective homomorphism of short exact sequences
\begin{equation} \label{eq:trop.of.unif.pts}
  \xymatrix @R=.2in{
    0 \ar[r] & {M'} \ar[r] \ar[d]^{\trop}_\cong
    & {E(K)} \ar[r] \ar[d]^{\trop}
    & {A(K)} \ar[r] \ar[d]^{\bar\tau} & 0 \\
    0 \ar[r] & {\trop(M')} \ar[r] 
    & {N_\Lambda} \ar[r] & {\Sigma_\Lambda} \ar[r] & 0
  }
\end{equation}
Both sequences~\eqref{eq:trop.of.unif} and~\eqref{eq:trop.of.unif.pts} are
compatible with extension of the ground field.

\paragraph[Uniformization and duality]
Let $T'$ be the split $K$-torus with character group $M'$.  
As explained in~\cite[\S6]{bosch_lutk:degenerating_avs},
the Raynaud cross~\eqref{eq:raynaud.cross} gives rise to a dual cross
\begin{equation} \label{eq:dual.raynaud.cross}
  \xymatrix @=.20in{
    & {T'^\an} \ar[d] & \\
    {M} \ar[r] & {E'^\an} \ar[r] \ar[d] & {\check A^\an} \\
    & {\check B^\an} &
  }
\end{equation}
where $\check A$ (resp.\ $\check B$) is the dual abelian variety of $A$
(resp.\ $B$), and $E'^\an$ is the universal cover of $\check A^\an$.
Note that
we have identifications $M = H_1(\check A^\an,\Z)$ and
$M' = H_1(A^\an,\Z)$, so duality interchanges character groups and
uniformizing lattices.  In the sequel we will need the following explicit
description of the composite homomorphism $T(K)\to E(K)\to A(K)$ in terms
of the data~\eqref{eq:dual.raynaud.cross}.  

We have $\Pic(\check A^\an) = \Pic(\check A)$ by analytic GAGA
theorems: see~\cite[\S1]{bosch_lutk:uniformizationII}.  Hence we have
a map on sheaf cohomology groups 
\[ H^1(\check A^\an,\underline K^\times)\To 
H^1(\check A^\an,\sO_{\check A^\an}^\times) = \Pic(\check A^\an) 
= \Pic(\check A) \]
induced by the inclusion of the constant sheaf
$\underline K^\times$ into $\sO_{\check A^\an}^\times$.
By~\cite[p.658]{bosch_lutk:degenerating_avs}, the line bundle
on $\check A$ corresponding to the image of a class in 
$H^1(\check A^\an, \underline K^\times)$ is translation-invariant, so we get a
homomorphism 
\begin{equation} \label{eq:pre.theta}
  H^1(\check A^\an,\underline K^\times) \To A(K).
\end{equation}
Identifying the sheaf cohomology group $H^1(\check A^\an,\underline K^\times)$ with
the singular cohomology group 
\[ H^1(\check A^\an,K^\times) = \Hom(H_1(\check A^\an,\Z),\,K^\times)
= \Hom(M, K^\times) = T(K) \]
and composing with~\eqref{eq:pre.theta} gives a homomorphism
\begin{equation} \label{eq:theta}
  \theta~:~ T(K) \To A(K). 
\end{equation}

\begin{lem} \label{lem:alternate.TtoA}
  The composition $T(K) \inject E(K) \to A(K)$ coincides with
  the homomorphism $\theta$ of~\eqref{eq:theta}.
\end{lem}

\pf Unwrapping the construction
of~\cite[\S4]{bosch_lutk:degenerating_avs}, one finds that the line bundle
on $\check A^\an$ corresponding to the image of a point 
$x\in T(K) = \Hom(M,K^\times)$ under the composition
$T(K)\to E(K)\to A(K)$ is obtained by
descending the trivial bundle on $E'^\an$ through the
exact sequence $0\to M\to E'^\an\to\check A^\an\to 0$
using the homomorphism $x: M \to K^\times$.  Regarding
$L = \theta(x)$ as a line bundle on $\check A^\an$, it is elementary to
check (as for topological line bundles)
that the canonical action of $M'$ on the pullback of $L$ to $E'$
is again given by the homomorphism $x$.\qed

\paragraph[Uniformization of Jacobians]
\label{par:unif.jacobian}
Let $X$ be a smooth, proper, geometrically connected $K$-curve
endowed with a split semistable $R$-model $\fX$.  Let $J$ be the Jacobian
of $X$.  Since $X$ has split semistable reduction, $J$ has split
semi-abelian reduction, so we have dual non-Archimedean uniformizations 
\begin{equation} \label{eq:jac.raynaud.cross}
  \xymatrix @=.20in{
    & {T^\an} \ar[d] & & &
    & {T'^\an} \ar[d] & \\
    {M'} \ar[r] & {E^\an} \ar[r] \ar[d] & {J^\an} & &
    {M} \ar[r] & {E'^\an} \ar[r] \ar[d] & {\check J^\an} \\
    & {B^\an} & & &
    & {\check B^\an} &
  }
\end{equation}
Recall that $M = H_1(\check J^\an,\Z)$ is the character lattice of $T$ and
$M' = H_1(J^\an,\Z)$ is the character lattice of $T'$.

\subparagraph \label{par:alternate.TtoJ}
Bosch and L\"utkebohmert~\cite[\S\S6--7]{bosch_lutk:uniformizationII} construct
the left square of~\eqref{eq:jac.raynaud.cross} directly in terms of the
curve $X$.  In their construction, the torus $T$ has character lattice
$H_1(X^\an,\Z)$, i.e.\ $M = H_1(X^\an,\Z)$.
The composition $\theta: T(K) \to E(K)\to J(K)$ can be described in a
manner analogous to Lemma~\ref{lem:alternate.TtoA} as the composition 
\begin{equation} \label{eq:theta.jac}
  T(K) = \Hom(H_1(X^\an,\Z),K^\times) = H^1(X^\an, K^\times)
  \to H^1(X^\an, \sO_{X^\an}^\times) = \Pic(X^\an), 
\end{equation}
which factors through $J(K)\subset\Pic(X) = \Pic(X^\an)$.
See Proposition~1.5 and Theorems~6.2 and~7.5 in \textit{loc.\ cit.\ }for
details.

\paragraph[The Abel-Jacobi map]
\label{par:compat.aj}
Suppose that $X(K)\neq\emptyset$, and let $\alpha:X\to J$ be an
Abel-Jacobi map.  The morphism
$\alpha^*:\check J\to J$ induced by Picard functoriality 
is the (inverse of the) canonical principal polarization of $J$.  
If $X$ has no $K$-rational points then $\alpha^*$ is defined by descent
theory after passing to a finite extension of $K$.  In either case,
$\alpha^*$ lifts in a unique way to an isomorphism
$\td\alpha^*: E'^\an\isom E^\an$, which then restricts to an isomorphism 
$\td\alpha^*: M\isom M'$.  In other words, we have an isomorphism of short
exact sequences
\begin{equation}  \label{eq:unif.polarization}
  \xymatrix @R=.2in{
    0 \ar[r] & M \ar[r] \ar[d]_\cong^{\td\alpha^*} 
    & {E'^\an} \ar[r] \ar[d]_\cong^{\td\alpha^*} 
    & {\check J^\an} \ar[r] \ar[d]_\cong^{\alpha^*} & 0 \\
    0 \ar[r] & {M'} \ar[r]
    & {E^\an} \ar[r] 
    & {J^\an} \ar[r] & 0 
  }
\end{equation}

\smallskip
We have the following Berkovich-analytic version of the classical
consequence of Abel's theorem which says that the Abel-Jacobi map realizes
the isomorphism $H_1(X(\C),\Z)\isom H_1(J(\C),\Z)$ in the complex setting.

\begin{prop} \label{prop:polariz.from.aj}
  Suppose that $X(K)\neq\emptyset$, and let 
  $\alpha: X\to J$ be an Abel-Jacobi map.  The homomorphism 
  on singular homology groups
  \[ \alpha_* ~:~ H_1(X^\an,\Z) \To H_1(J^\an,\Z) = M' \]
  coincides with the isomorphism
  $\td\alpha^*: M\isom M'$ of~\parref{par:compat.aj}
  under the identification $M = H_1(X^\an,\Z)$
  of~\parref{par:alternate.TtoJ}.  In particular, $\alpha_*$ is an
  isomorphism. 
\end{prop}

\pf The polarization $\alpha:\check J\isom J$ also induces an isomorphism
of Raynaud extensions
\begin{equation} \label{eq:raynaud.extn.polarization} \xymatrix @R=.2in{
  0 \ar[r] & {(T')^\an} \ar[d]^{\beta^*}_\cong \ar[r] 
  & {E'} \ar[d]^{\td\alpha^*}_\cong \ar[r] 
  & {\check B^\an} \ar[d]_\cong \ar[r] & 0 \\
  0 \ar[r] & {T^\an} \ar[r] & {E} \ar[r] & {B^\an} \ar[r] & 0 
}\end{equation}
Since polarizations are symmetric, 
by~\cite[Proposition~6.10(c)]{bosch_lutk:degenerating_avs} the map on
character groups $\beta: M\isom M'$ induced by the map
$\beta^*:T'\isom T$ coincides with $\td\alpha^*$
(note $\beta^*$ is algebraic because analytic
characters are algebraic). Therefore we may equivalently prove that 
this homomorphism $\beta$ on character groups coincides with $\alpha_*$. 

It is straightforward to show that the square
\[\xymatrix @=.2in{
  {T'(K)} \ar[r]^\theta \ar[d]_{\gamma} & 
  {\check J(K)} \ar[d]^{\alpha^*} \\
  {T(K)} \ar[r]^\theta & {J(K)} 
}\]
is commutative, where $\gamma:T'\to T$ is the homomorphism induced by
the map $\alpha_*:H_1(X^\an,\Z)\to H_1(J^\an,\Z)$ on characters, and the
maps $\theta$ are defined in~\eqref{eq:theta} 
and~\eqref{eq:theta.jac}.
On the other hand, the diagram
\[\xymatrix @=.2in{
  {T'(K)} \ar[r] \ar[d]_{\beta^*} &
  {E'(K)} \ar[r] \ar[d]^{\td\alpha^*} &
  {\check J(K)} \ar[d]^{\alpha^*} \\
  {T(K)} \ar[r] & {E(K)} \ar[r] & {J(K)}
}\]
made from~\eqref{eq:unif.polarization} and~\eqref{eq:raynaud.extn.polarization}
is also commutative, and the top (resp.\ bottom) row is again $\theta$ by
Lemma~\ref{lem:alternate.TtoA} (resp.~\parref{par:alternate.TtoJ}).
It follows that $\phi\coloneq\gamma-\beta^*: T'\to T$ takes $T'(K)$ into
the discrete subgroup $\ker(\theta) = T(K)\cap M'$.  Since $M'$ is a lattice
in $E$, we have
$\phi(T(K))\cap T_0(K) = \{0\}$, where $T_0\subset T^\an$ is the affinoid
torus $\trop\inv(0)$.  It is an
elementary exercise to show that any algebraic homomorphism of tori 
$\phi: T'\to T$ such that $\phi(T'(K))\cap T_0(K) = \{0\}$ is trivial, so
$\beta^* = \gamma$ and hence $\beta = \alpha_*$.\qed

\section{Retraction of divisors and divisor classes}
\label{sec:retraction.compat}

Recall that $K$ denotes a complete non-Archimedean field.

\paragraph[Curves and their skeleta]
At this point it is necessary to recall some of the analytic structure of 
$K$-curves.  The original reference is Berkovich's
book~\cite[Chapter~4]{berkovich:analytic_geometry}.
Our main reference is~\cite{bpr:analytic_curves}; see
also~\cite{temkin:berkovich_spaces}, Thuillier's
thesis~\cite{thuillier:thesis}, and  
Ducros' forthcoming book~\cite{ducros:curves_book}.
Let $X$ be a smooth, proper, geometrically connected
$K$-curve endowed with a split semistable $R$-model $\fX$.
This means that the special fiber $\fX_k$ of our model $\fX$ is
geometrically reduced, geometrically connected, and 
generically smooth, and for any singular point $x\in\fX_k$, the completed
local ring $\hat\sO_{\fX,x}$ is isomorphic to $R\ps{s,t}/(st-\varpi_x)$ for
some nonzero element $\varpi_x\in R$ with $\val(\varpi_x) > 0$.  The
\emph{incidence graph} $\Gamma = \Gamma_\fX$ associated to $\fX$ is the
metric $\Lambda$-graph (together with an associated $\Lambda$-rational model $G$) constructed as follows.  
The vertices of $\Gamma$ (or more properly, of $G$) are the irreducible components of $\fX_k$, and for each node $x\in\fX_k$,
there is an edge $e_x$ of length $\ell(e_x) \coloneq \val(\varpi_x)$ in
$\Gamma$ connecting the irreducible components containing $x$ (this is a
loop edge if $x$ is contained in a single irreducible component).  There
is a canonical inclusion $\Gamma\inject X^\an$ and a deformation
retraction $\tau: X^\an\to\Gamma$.  We identify $\Gamma$ with its image in
$X^\an$, called the \emph{skeleton} of $X$ (or of $X^\an$) associated to
$\fX$.  We have $\tau(X(K))\subset\Gamma(\Lambda)$, and 
$\tau: X(K)\to\Gamma(\Lambda)$ is surjective when $k$ is
algebraically closed (or more generally when every irreducible component
of $\fX_k$ contains a smooth $k$-rational point;
see~\cite[Proposition~10.1.40(b)]{liu:algebraic_geometry}).
Let $\Div_K(X)$ be the group of divisors of $X$ supported on $X(K)$.
Extending $\tau$ by linearity gives a homomorphism
\[ \tau_* ~:~ \Div_K(X) \To \Div_\Lambda(\Gamma). \]

The above constructions are compatible with extension of the ground field,
in the following sense.  Let $K'$ be a complete valued field extension of
$K$ and let $R'$ be its valuation ring.  Then $\fX_{R'}$ is a split
semistable $R'$-model of $X_{K'}$.  One checks that the associated
skeleton $\Gamma'\subset X_{K'}^\an$ maps isomorphically onto $\Gamma$
under the structure morphism $X_{K'}^\an\to X_K^\an$, and that the
following squares commute:
\[\xymatrix @=.2in{
  {\Gamma'} \ar[r] \ar[d]^\cong
  & {X_{K'}^\an} \ar[d] 
  & & {X_{K'}^\an} \ar[r]^{\tau} \ar[d]
  & {\Gamma'} \ar[d]^\cong \\
  {\Gamma} \ar[r] & {X^\an} & & {X^\an} \ar[r]^\tau & {\Gamma}
}\]

\paragraph
Let $J$ be the Jacobian of $X$.  It has a uniformization 
\[ \xymatrix @=.20in{
  & {T^\an} \ar[d] & \\
  {M'} \ar[r] & {E^\an} \ar[r]^p \ar[d] & {J^\an} \\
  & {B^\an} &
}\]
as in~\parref{par:unif.jacobian}, where $T$ has character lattice
$M = H_1(X^\an,\Z) = H_1(\Gamma,\Z)$.  
Likewise the Jacobian of $\Gamma$ has a uniformization of the form
$0\to M\to N_\R\overset{q}\to\Jac(\Gamma)\to 0$, where 
$N = \Hom(M,\Z) = H^1(\Gamma,\Z)$ as usual; see
Proposition~\ref{prop:div.unif.jac2}. 
The goal of this section is to define a \emph{specialization map} 
$\bar\tau\p: J^\an\to\Jac(\Gamma)$
which is compatible with the complementary descriptions of $J$ and of
$\Jac(\Gamma)$ in terms of divisor classes and uniformizations.  This will
be a key ingredient in the proof of
Theorem~\ref{thm:skeleta.compatibility}, which says that 
$\Jac(\Gamma)\cong\Sigma(J)$ canonically; under this identification, we
will have $\bar\tau\p=\bar\tau$ (see
Corollary~\ref{cor:skeleta.compatibility}).

\begin{prop} \label{prop:define.taubar}
  There exists a unique surjective homomorphism 
  $\bar\tau\p: J(K)\surject\Jac_\Lambda(\Gamma)$ making the following two
  squares commute:
  \begin{equation} \label{eq:define.taubar} \xymatrix @=.25in{
    {\Div_K^0(X)} \ar[r] \ar[d]_{\tau_*} & {J(K)} \ar[d]^{\bar\tau\p} & &
    {E(K)} \ar[r]^p \ar[d]_{\trop} & {J(K)} \ar[d]^{\bar\tau\p} \\
    {\Div_\Lambda^0(\Gamma)} \ar[r] & {\Jac_\Lambda(\Gamma)} & &
    {N_\Lambda} \ar[r]_(.4)q & {\Jac_\Lambda(\Gamma)}
  }\end{equation}
\end{prop}

\pf Suppose first that $K$ is algebraically closed, so that
$\Div_K(X) = \Div(X)$.
Let $\Prin(X) = \ker(\Div^0(X)\to J(K))$ denote the group of
principal divisors.  We have
$\tau_*(\Prin(X))\subset\Prin_\Lambda(\Gamma)$ by the Poincar\'e-Lelong
formula (see Remark~\ref{rem:poincare.lelong} below), so there exists a unique map
$\bar\tau\p:J(K)\to\Jac_\Lambda(\Gamma)$  
making the left square in~\eqref{eq:define.taubar} commute.
We claim that $\bar\tau\p\circ p = q\circ\trop$, so that the right square
commutes as well.
We have $E(K) = T(K) + J_0(K)$ (where $J_0 = \bar\tau\inv(0)$
as in~\parref{par:unif.basics.2}), so it suffices to show
$\bar\tau\p(p(x)) = q(\trop(x))$ for $x\in T(K)$ and $x\in J_0(K)$.

\subparagraph Let $x\in J_0(K)\subset J(K)$ and let $\sL$ be a line bundle
representing $x$.  Then $\sL$ is represented by a divisor $D$ such that
$\tau_*(D) = 0$ by~\cite[Theorem~5.1(c)]{bosch_lutk:uniformizationII}, so 
$\bar\tau\p(x) = 0$.  By definition we have $\trop(x) = 0$ as well.

\subparagraph Choose a $\Lambda$-rational model $G$ of $\Gamma$ with no
loop edges~\parref{par:graph.basics},  and let $\cV(G)$ (resp.\ $\cE(G)$)
denote the set of vertices (resp.\ edges) of $G$.  For $v\in \cV(G)$ let 
$U_v$ be the open star around $v$ in $\Gamma$, i.e.\ 
the union of $\{v\}$ and all open edges adjacent to $v$.
Let $\cU$ denote the open cover $\{U_v~:~v\in \cV(G)\}$ of $\Gamma$.  
Let $U_v' = \tau\inv(U_v)$; this is an open analytic domain in $X^\an$,
and $\cU' = \{U_v'~:~v\in \cV(G)\}$ is an 
open cover.  Choose orientations of the edges of $G$,
so that $G$ becomes a simplicial complex.
Since each $U_v$ and $U_v'$ is contractible, we have
canonical isomorphisms 
\[ \check H^1(\cU',K^\times) = H^1(X^\an, K^\times) \sptxt{and}
\check H^1(\cU,\Lambda) = H^1(\Gamma,\Lambda). \]
Let $H\in\check C^1(\cU',K^\times)$ be a \v Cech $1$-cocycle, so $H$
amounts to a function $H: \cE(G)\to K^\times$
and defines a class $[H]\in H^1(X^\an,K^\times) = \Hom(M,K^\times) =
T(K)$.  Extending $H$ linearly defines a homomorphism 
$H: C_1(\cV(G),\Z)\to K^\times$ on simplicial $1$-chains; the
restriction of this homomorphism to $H_1(\Gamma,\Z)$ is the element of
$\Hom(H_1(\Gamma,\Z),K^\times) = H^1(\Gamma,K^\times) = H^1(X^\an,K^\times)$ 
corresponding to the class of $H$.
Let $h = \val\circ H: \cE(G)\to\Lambda$.
This is a cocycle in $\check C^1(\cU,\Lambda)$, and
we have an equality of cohomology classes 
$[h] = \trop([H])\in H^1(\Gamma,\Lambda) = N_\Lambda$. 

Let $\sL$ be the algebraization of the analytic line bundle on $X$ given
by the cocycle 
$H\in\check C^1(\cU',K^\times)\subset\check C^1(\cU',\sO_{X^\an}^\times)$.
In other words, $\sL$ is the image of $[H]$ under the composition
\[ H^1(X^\an,K^\times) \to H^1(X^\an,\sO_{X^\an}^\times) = \Pic(X^\an)
 = \Pic(X). \]
Regarding $[H]$ as an element of $T(K)$,
by~\parref{par:alternate.TtoJ} we have that $p([H])$ is the isomorphism 
class of $\sL$.  Since $\sL^\an$ was
constructed using the \v Cech cocycle $H$, we have distinguished trivializations 
$\sO_{\sL^\an}|_{U_v'} = \sO_{U_v'}$ for each $v\in \cV(G)$, where 
$\sO_{\sL^\an}$ is the sheaf of sections of $\sL^\an$; the cocycle $H$ determines
the transition functions on $U_v'\cap U_w'$.
Similarly, let $L$ be the line bundle
on $\Gamma$ given by the cocycle 
$h\in\check C^1(\cU,\Lambda)\subset\check C^1(\cU,\sO_\Lambda)$, so 
$q([h])$ is the isomorphism class of $L$; see Proposition~\ref{prop:pic.from.H1}.
We have distinguished trivializations 
$\sO_{L,\Lambda}|_{U_v} = \sO_\Lambda|_{U_v}$ for each $v\in \cV(G)$, and
$h$ determines the transition functions.

Let $f$ be a nonzero meromorphic section of $\sL$ and let 
$D = \div(f)$, so $\sL = p([H]) = [D]$.  We may identify $f$ with a
collection $\{f_v~:~v\in \cV(G)\}$, where $f_v$ is a meromorphic
function on $U_v'$ and for each edge $e = \overrightarrow{vw}$
we have $H(e)\,f_v = f_w$ on $\tau\inv(e^\circ)$.  Let 
$F_v$ be the restriction of $-\log |f_v|$ to $\Gamma$.  By the
Poincar\'e-Lelong formula (see Remark~\ref{rem:poincare.lelong} below),
$F_v$ is a tropical meromorphic function on $U_v$, and 
$\div(F_v) = \tau_*(\div(f_v))$.  For $e = \overrightarrow{vw}$
we have $h(e) + F_v = F_w$ on $e^\circ$, 
so the collection $\{F_v~:~v\in \cV(G)\}$
glues to give a meromorphic section $F$ of $L$.  By construction
we have $\tau_*(D) = \div(F)$.  
By definition $[L] = [\div(F)]$, so
\[ q(\trop([H])) = q([h]) = [L] = [\tau_*(D)] = \bar\tau\p([D]) =
\bar\tau\p(p([H])). \]
Therefore the right square in~\eqref{eq:define.taubar} commutes.

\subparagraph
Now suppose that $K$ is not necessarily algebraically closed.  Let $K'$
denote the completion of an algebraic closure of $K$ and let 
$\Lambda' = \val(K'^\times)$.  Define 
$\bar\tau\p: J(K)\to\Jac_{\Lambda'}(\Gamma)$ to be the composition
\[ J(K) \inject J(K') \overset{\bar\tau\p}\To \Pic_{\Lambda'}^0(\Gamma). \]
Then the following two diagrams commute:
\[\xymatrix @R=.2in{
  {\Div_K^0(X)} \ar[r] \ar[d] & {J(K)} \ar[d] & &
  {E(K)} \ar[r] \ar[d] & {J(K)} \ar[d] \\
  {\Div_{K'}^0(X)} \ar[r] \ar[d]_{\tau_*} & {J(K')} \ar[d]^{\bar\tau\p} & &
  {E(K')} \ar[r] \ar[d]_{\trop} & {J(K')} \ar[d]^{\bar\tau\p} \\
  {\Div_{\Lambda'}^0(\Gamma)} \ar[r] & {\Jac_{\Lambda'}(\Gamma)} & &
  {N_{\Lambda'}} \ar[r] & {\Jac_{\Lambda'}(\Gamma)}
}\]
Since $\trop(E(K)) = N_\Lambda$ and 
$N_\Lambda$ surjects onto $\Jac_\Lambda(\Gamma)$, the image of
$\bar\tau\p:J(K)\to\Jac_{\Lambda'}(\Gamma)$ is equal to the subgroup
$\Jac_\Lambda(\Gamma)$.  Therefore 
$\bar\tau\p:J(K)\surject\Jac_\Lambda(\Gamma)$ makes the
squares~\eqref{eq:define.taubar} commute.  The uniqueness of $\bar\tau\p$
follows from the surjectivity of $E(K)\to J(K)$.\qed

\begin{rem} \label{rem:poincare.lelong}
  In the above proof we used
  the non-Archimedean Poincar\'e-Lelong formula, which can be formulated as
  follows.  Let $X$ be a smooth curve and $\Gamma$ a skeleton of $X^\an$,
  as above.  Let $f$ be a nonzero meromorphic function on $X$ and let 
  $F = -\log |f|$ restricted to $\Gamma$.  Then $F$ is a
  $\Lambda$-rational meromorphic function on $\Gamma$ and $\div(F)$ is the
  retraction of $\div(f)$.  A statement and proof in our language can be found
  in~\cite[\S5]{bpr:analytic_curves}, where it is called the slope
  formula.  The result is originally due to
  Thuillier~\cite[Proposition~3.3.15]{thuillier:thesis}, where it is
  stated in the language of non-Archimedean potential theory.
  Thuillier's theorem is also proved in greater generality; in particular
  it applies to all smooth analytic curves $X$, not just algebraic ones.
\end{rem}

\paragraph \label{par:Mp.to.M}
The second square in~\eqref{eq:define.taubar} extends to a
homomorphism of short exact sequences
\begin{equation} \label{eq:trop.of.unif2}
  \xymatrix @R=.25in {
    0 \ar[r] & {M'} \ar[r] \ar[d] & {E(K)} \ar[r]^p \ar[d]_{\trop}
    & {J(K)} \ar[r] \ar[d]_{\bar\tau\p} & 0 \\
    0 \ar[r] & {M} \ar[r] & {N_\Lambda} \ar[r]_(.4)q 
    & {\Jac_\Lambda(\Gamma)} \ar[r] & 0
  }
\end{equation}
In particular $\trop(M')\subset M$, so there is a natural homomorphism
\begin{equation} \label{eq:skel.to.jac}
  \pi ~:~ \Sigma(J) = N_\R/\trop(M') \To N_\R/M = \Jac(\Gamma)
\end{equation}
such that $\pi\circ\bar\tau = \bar\tau\p$.
Composing with $N_\R\to\Sigma(J)$, we have proved:

\begin{cor} \label{cor:taubar.2}
  There is a unique surjective map $\bar\tau\p: J^\an\to\Jac(\Gamma)$
  making the square
  \begin{equation} \label{eq:taubar.2} \xymatrix @=.25in{
    {E^\an} \ar[r]^p \ar[d]_{\trop} & {J^\an} \ar[d]^{\bar\tau\p} \\
    {N_\R} \ar[r]_(.4)q & {\Jac(\Gamma)}
  }\end{equation}
  commute.  Restricting to $K$-points yields the second square
  in~\eqref{eq:define.taubar}. 
\end{cor}

\section{The skeleton of the Jacobian is the Jacobian of the skeleton}
\label{sec:skel.jacobian}

We continue to use the notations in \S\ref{sec:retraction.compat}.  In
particular, $X$ is a curve with split semistable reduction, $J$ is its
Jacobian, and $M = H_1(X^\an,\Z) = H_1(\check J^\an,\Z)$ is the character
lattice of the torus in the uniformization of $J$. 
Let $\Gamma\subset X^\an$ be the skeleton of $X$ associated to a split
semistable model $\fX$, and let $\tau:X^\an\to\Gamma$ be the retraction
map. 

The following compatibilty of Abel-Jacobi maps is the key ingredient in
the proof of Theorem~\ref{thm:skeleta.compatibility}.  Its proof uses
Proposition~\ref{prop:define.taubar} in an essential way.

\begin{prop} \label{prop:compat.aj}
  Suppose that $X(K)\neq\emptyset$.  Let $P\in X(K)$ and let 
  $p = \tau(P)\in\Gamma(\Lambda)$.  Let $\alpha_P: X\to J$ and 
  $\alpha_p: \Gamma\to\Jac(\Gamma)$ be the Abel-Jacobi maps based at $P$ and
  $p$, respectively.   Then the following square commutes: 
  \begin{equation} \label{eq:aj.commutes}\xymatrix @=.25in{
      {X^\an} \ar[r]^{\alpha_P} \ar[d]_\tau & 
      {J^\an} \ar[d]^{\bar\tau\p} \\
      {\Gamma} \ar[r]_(.4){\alpha_p} & {\Jac(\Gamma)}
    }\end{equation}
  where $\bar\tau\p$ is the map of Corollary~\ref{cor:taubar.2}.
\end{prop}

\pf Let $K'$ be a complete valued field extension of $K$.  Define
$\td\alpha_P: X(K')\to\Div_{K'}^0(X_{K'})$ by $Q\mapsto (Q)-(P)$ and
$\td\alpha_p: \Gamma\to\Div^0(\Gamma)$ by $q\mapsto (q)-(p)$.  
The commutativity of the left square in the diagram
\[\xymatrix @R=.2in{
  {X(K')} \ar[r]^(.43){\td\alpha_P} \ar[d]_\tau &
  {\Div^0_{K'}(X_{K'})} \ar[r] \ar[d]_{\tau_*} &
  {J(K')} \ar[d]^{\bar\tau\p} \\
  {\Gamma} \ar[r]_(.4){\td\alpha_p} &
  {\Div^0(\Gamma)} \ar[r] & {\Jac(\Gamma)}
}\]
is obvious, and the commutativity of the right square
follows from the commutativity of the left square
in~\eqref{eq:define.taubar} after extending the ground field.  The
compositions of the horizontal arrows in the above diagram
are the maps $\alpha_P$ and $\alpha_p$,
so the outer square and the top square of the diagram
\[\xymatrix @R=.2in{
  {X(K')} \ar[r]^{\alpha_P} \ar[d] & {J(K')} \ar[d] \\
  {X^\an} \ar[r]^{\alpha_P} \ar[d]_\tau & {J^\an} \ar[d]^{\bar\tau\p} \\
  {\Gamma} \ar[r]^(.4){\alpha_p} & {\Jac(\Gamma)}
}\]
are commutative.  It follows that the bottom square commutes as well,
since any point $x\in X^\an$ is the image of a morphism 
$\Spec(K')\to X^\an$ for some $K'$.\qed

\paragraph
If $X(K)\neq\emptyset$ we let $\alpha: X\to J$ be an Abel-Jacobi map,
and in general we let $\alpha^*:\check J\isom J$ be the canonical principal
polarization, as in~\parref{par:compat.aj}.  Let
$\td\alpha^*: E'^\an\isom E^\an$ be the lift of $\alpha^*$ to universal
covers; this restricts to an isomorphism
$\td\alpha^*: M\isom M' = H_1(J^\an,\Z)$ on homology groups.  Letting
$N = \Hom(M,\Z)$, the composition $\trop\circ\td\alpha^*:M\to N_\Lambda$
defines a pairing 
\begin{equation} \label{eq:analytic.monodromy}
  \angles{\scdot,\scdot}_\an ~:~ M\times M \To \Lambda,
\end{equation}
which we call the \emph{analytic monodromy pairing}.

Since $M = H_1(X^\an,\Z) = H_1(\Gamma,\Z)$, the Jacobian of
$\Gamma$ has a uniformization of the form 
\[ 0 \To M \To N_\R \To \Jac(\Gamma) \To 0, \]
where the homomorphism $M\to N_\R$ comes from the edge length pairing
$\angles{\scdot,\scdot}:M\times M\to\Lambda$
of~\eqref{eq:edge.length.pairing}. 

\begin{rem} \label{rem:analytic.groth.monodromy}
  The analytic monodromy pairing $\angles{\scdot,\scdot}_\an$ coincides
  with Grothendieck's monodromy 
  pairing in the case $\Lambda = \Z$.  This is proved
  in~\cite{coleman:monodromy_pairing}.
  Note that this fact is unrelated to the more difficult
  Theorem~\ref{thm:analytic.monodromy.pairing}, which says that the edge length
  pairing~\eqref{eq:edge.length.pairing.1} coincides with the monodromy
  pairing.  Indeed, the former fact holds for arbitrary abelian varieties,
  and the standard proofs of the latter make no reference to
  the non-Archimedean uniformization theory of $J$.
\end{rem}

\smallskip
We are now able to prove Theorem~\ref{thm:skeleta.compatibility}:

\begin{thm*}[Theorem~\ref{thm:skeleta.compatibility}]
  The analytic monodromy pairing
  $\angles{\scdot,\scdot}_\an$ defined in~\parref{par:compat.aj} coincides
  with the edge length pairing $\angles{\scdot,\scdot}$
  under the identification $M = H_1(X^\an,\Z) = H_1(\Gamma,\Z)$
  of~\parref{par:alternate.TtoJ}.
\end{thm*}

\pf We observed in~\parref{par:Mp.to.M} that 
$\trop: E(K)\to N_\Lambda$ takes $M'$ into $M = H_1(\Gamma,\Z)$.  
The Theorem is equivalent to the assertion that the composition
$\trop\circ\td\alpha^*: M\to M$ is the identity, which
can be checked after extending the ground field.  Hence we may
assume that $X$ has a $K$-rational point $P$.  Let
$p = \tau(P)\in\Gamma(\Lambda)$ and let
$\alpha_P: X\to J$ (resp.\ $\alpha_p:\Gamma\to\Jac(\Gamma)$)
be the Abel-Jacobi map based at $P$ (resp.\ at $p$).
Applying $H_1(\scdot,\Z)$ to~\eqref{eq:aj.commutes}, we get the
commutativity of the left square in
\[\xymatrix @R=.2in{
  {H_1(X^\an,\Z)} \ar[r]^{(\alpha_P)_*} \ar[d]_{\tau_*}^\cong 
  & {H_1(J^\an,\Z)} \ar[d]_{\bar\tau\p_*} \ar[r]
  & {E^\an} \ar[d]^{\trop} \\
  {H_1(\Gamma,\Z)} \ar[r]_(.45){(\alpha_p)_*} 
  & {H_1(\Jac(\Gamma),\Z)}  \ar[r]
  & {N_\R}
}\]
The right square commutes by basic point-set topology: note that the map 
$\trop:E^\an\to N_\R$ is the (unique) lift of 
$\bar\tau\p:J^\an\to\Jac(\Gamma)$ to universal covers.

The map $(\alpha_P)_*$ is the isomorphism
$\td\alpha^*: M\to M' = H_1(J^\an,\Z)$ by
Proposition~\ref{prop:polariz.from.aj}.  One checks that the composition
of the bottom horizontal arrows is the homomorphism $M\to N_\R$ defined by
the edge length pairing: see~\parref{par:tropical.aj} and
Lemma~\ref{lem:graph.aj.kpi1}. Hence
$\trop\circ\td\alpha^*$ takes $H_1(X^\an,\Z)$ into $H_1(\Gamma,\Z)$ via
the isomorphism $\tau_*$ functorially induced by the retraction $\tau$.
We are using $\tau_*$ to identify $H_1(\Gamma,\Z)$ with $M =
H_1(X^\an,\Z)$, so $\trop\circ\td\alpha^*$ is the identity map under our
identifications.\qed 

\begin{rem} \label{rem:second.proof}
  Suppose that $\Lambda = \Z$.  
  Theorem~\ref{thm:skeleta.compatibility} says that the edge length
  pairing $\angles{\scdot,\scdot}$ coincides with the analytic monodromy
  pairing $\angles{\scdot,\scdot}_\an$, which is the same as
  Grothendieck's monodromy pairing by
  Remark~\ref{rem:analytic.groth.monodromy}.  We thus have a new proof of
  Theorem~\ref{thm:analytic.monodromy.pairing}, which is very different
  from the standard ones as it does not use arithmetic intersection theory
  or the Picard-Lefschetz formula.
\end{rem}

\paragraph
We can rephrase Theorem~\ref{thm:skeleta.compatibility} by saying that the
uniformizations of $J$ and of $\Jac(\Gamma)$ are ``compatible under
tropicalization'' in the sense that the following diagram commutes:
\[\xymatrix @R=.25in {
    0 \ar[r] & {M} \ar[r] \ar@{=}[d] & {E^\an} \ar[r] \ar[d]_{\trop}
    & {J^\an} \ar[r] \ar[d]_{\bar\tau\p} & 0 \\
    0 \ar[r] & {M} \ar[r] & {N_\R} \ar[r]
    & {\Jac(\Gamma)} \ar[r] & 0
}\]
Here we have used $\td\alpha^*$ to identify $M' = H_1(J^\an,\Z)$
with $M = H_1(X^\an,\Z)=H_1(\Gamma,\Z)$.  Since 
$\Sigma(J) = N_\R/\trop(M)$, the map 
$\pi: \Sigma(J)\to\Jac(\Gamma)$ of~\eqref{eq:skel.to.jac} is an
isomorphism, so we have proved:

\begin{cor} \label{cor:skeleta.compatibility}
  There is a canonical isomorphism of the skeleton $\Sigma(J)$ of the
  Jacobian of $X$ with the Jacobian $\Jac(\Gamma)$ of the skeleton of
  $X$ as principally polarized real tori, with the polarizations defined
  by the analytic monodromy pairing and the edge length pairing,
  respectively.  The map $\bar\tau\p:J^\an\to\Jac(\Gamma)$ of 
  Corollary~\ref{cor:taubar.2} coincides with the retraction to the
  skeleton $\bar\tau: J^\an\to\Sigma(J)$ under this identification.
\end{cor}

\begin{rem} \label{rem:period.matrix}
  Suppose that $X$ is a Mumford curve, so that its Jacobian $J$ has
  totally degenerate reduction.  In this case $E = T$ and the map
  $M\to T(K) = \Hom(M,K^\times)$ is given by a pairing
  $[\scdot,\scdot]:M\times M\to K^\times$.
  Choosing a basis for $M$, one can represent this pairing by a matrix
  $m$ with entries in $K^\times$.  We call $m$ the \emph{period matrix} of
  $J$.  Similarly, we call the matrix for the pairing
  $\angles{\scdot,\scdot}$ the \emph{period matrix} of the 
  principally polarized tropical abelian variety $\Jac(\Gamma)$.
We have $\angles{\scdot,\scdot}_\an = \val\circ [\scdot,\scdot]$
  by definition, so Theorem~\ref{thm:skeleta.compatibility} says that the
  period matrix of $\Jac(\Gamma)$ is the valuation of the period matrix of
  $J$.  Cf.\ Remark~\ref{rem:trop.of.aj}.
\end{rem}

\begin{cor} \label{cor:ker.taubarp}
  The inverse image of zero under the specialization map 
  $\bar\tau\p: J^\an\to\Jac(\Gamma)$ is equal to the analytic domain 
  $J_0\subset J^\an$ of~\parref{par:unif.basics.2}.
\end{cor}

\pf By definition $J_0 = \bar\tau\inv(0)$, so this follows from 
Corollary~\ref{cor:skeleta.compatibility}.\qed

\section{Surjectivity of principal divisors}

In this section we assume that $K$ is an \emph{algebraically closed}
complete non-Archimedean field.  
Otherwise we keep the notations of
\S\S\ref{sec:retraction.compat}--\ref{sec:skel.jacobian}.

\paragraph
Since $K$ is algebraically closed we have
$\Div_K^0(X) = \Div^0(X)$, which surjects onto $J(K)$ and onto
$\Div^0_\Lambda(X)$; taking kernels of the
left square in~\eqref{eq:define.taubar}, we get a homomorphism of short
exact sequences
\begin{equation} \label{eq:div.jac.seq}
  \xymatrix @R=.25in{
    0 \ar[r] & {\Prin(X)} \ar[r] \ar[d]^{\tau_*}
    & {\Div^0(X)}  \ar[r] \ar[d]^{\tau_*} 
    & {J(K)} \ar[r] \ar[d]^{\bar\tau\p} & 0 \\
    0 \ar[r] & {\Prin_\Lambda(\Gamma)} \ar[r] 
    & {\Div^0_\Lambda(\Gamma)} \ar[r]
    & {\Jac_\Lambda(\Gamma)} \ar[r] & 0
  }  
\end{equation}
with surjective middle and right vertical arrows.
By the Poincar\'e-Lelong formula (Remark~\ref{rem:poincare.lelong}), if
$f$ is a nonzero meromorphic function on $X$ and
$F = -\log|f|\big|_\Gamma$ then $F$ is a $\Lambda$-rational meromorphic
function on $\Gamma$ and
$\tau_*(\div(f)) = \div(F)$; this explains
the left vertical arrow.
We can now prove the first part of Theorem~\ref{thm:prin.surj}.

\begin{thm*}[Theorem~\ref{thm:prin.surj}]
  The map on principal divisors $\tau_*: \Prin(X)\to\Prin_\Lambda(\Gamma)$
  is surjective.
\end{thm*}

\pf Let 
$\Div^{(0)}(X) = \ker(\tau_*:\Div^0(X)\to\Div^0_\Lambda(\Gamma))$.  By
Corollary~\ref{cor:ker.taubarp} we have $J_0(K) = \ker(\bar\tau\p)$.
Applying the snake lemma to the sequence~\eqref{eq:div.jac.seq}, it is
enough to prove that the map $\Div^{(0)}(X)\to J_0(K)$ is surjective.  But
this is exactly~\cite[Theorem~5.1(c)]{bosch_lutk:uniformizationII}.\qed

\begin{cor*}[Corollary~\ref{cor:merom.surj}]
  Let $F:\Gamma\to\R$ be a continuous function.  Then there exists a
  nonzero rational function $f\in K(X)$ such that 
  $F = -\log|f|\,\big|_\Gamma$ if and only if $F$ is a $\Lambda$-rational tropical meromorphic function.
\end{cor*}

\pf We have already remarked that the restriction of $-\log|f|$ to
$\Gamma$ is meromorphic and $\Lambda$-rational.  Conversely, let
$F$ be a $\Lambda$-rational tropical meromorphic function on $\Gamma$ and 
let $D = \div(F)$.  By Theorem~\ref{thm:prin.surj} there exists a
nonzero meromorphic function $f$ on $X$ such that 
$\tau_*(\div(f)) = D$.  Letting $F' = -\log|f|\big|_\Gamma$, we have
$\div(F') = D = \div(F)$, so $F-F' = \lambda\in \Lambda$ is a constant
function.  Choosing $a\in K^\times$ with valuation $\lambda$, we have
$F = -\log|af|\big|_\Gamma$.\qed

\paragraph
Let $\bar D = \bar D_1 - \bar D_2\in\Prin_\Lambda(\Gamma)$, where
$\bar D_1$ and $\bar D_2$ are effective divisors of degree $n\geq 0$.  The 
surjectivity of $\tau_*:\Prin(X)\to\Prin_\Lambda(\Gamma)$ does \emph{not}
imply that there 
exist effective divisors $D_1,D_2$ on $X$ of degree $n$ such that $D = D_1-D_2$ is
principal and $\tau_*(D_i) = \bar D_i$ for $i=1,2$ (see
Remark~\ref{rem:naive.false} for a counterexample).
This is because the map $\tau_*:\Div^0(X)\to\Div^0_\Lambda(\Gamma)$ allows for 
cancellation: for example, if $x,y\in X(K)$ retract to the
same point of $\Gamma$ then $\tau_*((x)-(y)) = 0$.  
The second part of Theorem~\ref{thm:prin.surj} states that we can limit
the amount of cancellation that happens:

\begin{thm*}[Theorem~\ref{thm:prin.surj}, continued]
  Suppose that $X$ has genus $g$.
  Let $\bar D = \bar D_1 - \bar D_2\in\Prin_\Lambda(\Gamma)$, where
  $\bar D_1,\bar D_2\in\Div^n_\Lambda(\Gamma)$ are effective divisors of
  degree $n$.  There exist effective divisors $D_1,D_2\in\Div^{n+g}(X)$
  of degree $n+g$ such that $D = D_1-D_2\in\Prin(X)$ and $\tau_*(D) = \bar D$.
\end{thm*}

\pf We will prove the following more precise statement.  
Let $D'\in\Div^0(X)$ be any divisor such that $\tau_*(D') = \bar D$.  We
claim that there exist effective divisors $E_1,E_2\in\Div(X)$ such that:
\begin{enumerate}
\item $\deg(E_1) = \deg(E_2) \leq g$,
\item $\tau_*(E_1 - E_2) = 0$, and
\item the divisor $D\coloneq D' + E_1 - E_2\in\Div^0(X)$ is
  principal.
\end{enumerate}
Chasing $D'$ around the diagram~\eqref{eq:div.jac.seq}
and using Corollary~\ref{cor:ker.taubarp}, one sees that the claim
amounts to proving that any element of  
$J_0(K)$ has the form $[E_1-E_2]$ for $E_1,E_2\in\Div(X)$
satisfying~(1) and~(2).

Let $\fX$ be a semistable model of $X$ whose special fiber has smooth
irreducible components.  As explained in~\parref{par:unif.basics}
and~\cite[\S5]{bosch_lutk:uniformizationII}, the $K$-analytic group
$J_0$ has semi-abelian reduction; in fact its special
fiber identified with the generalized Jacobian $\Jac(\fX_k)$ of the special fiber
$\fX_k$.  Hence there is a reduction map $\red: J_0(K)\to\Jac(\fX_k)(k)$.
Let $\bar C_1,\ldots,\bar C_r$ be the irreducible components of $\fX_k$.  
There is a canonical short exact sequence
\[ 0 \To \bar T \To \Jac(\fX_k) \overset{\rho}\To \prod_{i=1}^r \Jac(\bar C_i) \To 0, \]
where $\bar T$ is the split $k$-torus with character lattice $M$; its 
dimension is equal to the cyclomatic number 
$t\coloneq\dim_\Q H_1(\Gamma,\Q)$.

Let $e_1,\ldots,e_t$ be a set of edges for some $\Lambda$-rational model of $\Gamma$ whose complement is a
maximal spanning tree.  
For $i=1,\ldots,t$ choose a $\Lambda$-point
$\epsilon_i$ in the interior of $e_i$, and set 
$R = \prod_{i=1}^t\tau\inv(\epsilon_i)$.  Then
$R$ is a product of $t$ affinoid annuli of modulus zero,
so its canonical reduction $\bar R$ is isomorphic to $\bG_{m,k}^t$.  
For $i=1,\ldots,r$ let $\zeta_i\in X^\an$ be the unique  point reducing to
the generic point of $\bar C_i$ and let $C_i = \tau\inv(\zeta_i)$, so 
$C_i$ is the set of all points of $X^\an$ reducing to smooth points of
$\bar C_i$.  As in~\cite[Situation~3.4$''$]{bosch_lutk:uniformizationII} we define
\[ R \times C^{(*)} \coloneq R\times \prod_{i=1}^r C_i^{(g_i)} \subset (X^\an)^{(g)}, \]
where $g_i$ is the genus of $\bar C_i$ and the superscripts denote
symmetric powers.  This is an affinoid space with canonical reduction 
$(R \times C^{(*)})^-=\bG_{m,k}^t\times\prod_{i=1}^r (\bar C_i{}^\sm)^{(g_i)}$.  
Points of $R\times C^{(*)}(K)$ can be interpreted as degree-$g$ effective divisors
on $X$, all of which have the same retraction
$\sum_{i=1}^t (\epsilon_i) + \sum_{j=1}^r g_i\,(\zeta_i)$
to $\Gamma$.  For a divisor $E_2\in R\times C^{(*)}(K)$
let $\phi_{E_2}: R\times C^{(*)}\to J^\an$ be the morphism defined on
$K$-points by $E_1\mapsto[E_1-E_2]$.  
The image of $\phi_{E_2}$ is contained in $J_0 = \ker(\bar\tau\p)$
by~\eqref{eq:div.jac.seq} and Corollary~\ref{cor:ker.taubarp}
since $E_1-E_2\in\ker(\tau_*)$ for all $E_1\in R\times C^{(*)}(K)$.
We will prove that for a suitable divisor $E_2\in R\times C^{(*)}(K)$, the map
$\phi_{E_2}: R\times C^{(*)}(K)\to J_0(K)$ is surjective. 

Let $E_2\in R\times C^{(*)}(K)$.  Define a map
$\bar\phi_{E_2}:(R \times C^{(*)})^-(k)\to\Jac(\fX_k)(k)$ as follows.
Let $\fX'$ be the semistable model of $X$ obtained by blowing
up the nodes of $\fX_k$ in such a way that each $\epsilon_i$ reduces to a
generic point of $\fX'_k$.  Then the points of  $(R\times C^{(*)})^-(k)$
can be interpreted as degree-$g$ divisors on  
$\fX'_k$ with support on the smooth locus.  Letting
$\bar E_2\in\Div^g(\fX'_k)$ be the reduction of $E_2$, we set
$\bar\phi_{E_2}(\bar E_1) = [\bar E_2 - \bar E_1]\in\Jac(\fX'_k)
= \Jac(\fX_k)$.  By~\cite[Theorem~5.1(c)]{bosch_lutk:uniformizationII} we
have $\bar\phi_{E_2}\circ\red = \red\circ \phi_{E_2}$.

Write $\bar E_2 = (\bar E\p_2,\bar E{}''_2)\in(R \times C^{(*)})^-(k)$,
where $\bar E\p_2\in\bar R(k)=\G_{m,k}^t(k)$ and 
$\bar E{}''_2\in\prod_{i=1}^r(\bar C_i{}^\sm)^{(g_i)}(k)$.
The map $\bar\phi_{E_2}$ identifies 
$\bG_{m,k}^t(k) = \bar R(k)\times\{\bar E{}''_2\}\subset(R \times C^{(*)})^-(k)$
with the toric part $\bar T(k)$ of $\Jac(\fX_k)$: see the end of \S4
in~\cite{bosch_lutk:uniformizationII}.  Hence
$\bar T(k)$ is contained in the image of $\bar\phi_{E_2}$.  
For a generic choice of $\bar E_2$, the subset
$\bar\phi_{E_2}(\{\bar E\p_2\}\times \prod_{i=1}^r (\bar C_i{}^\sm)^{(g_i)})\subset\Jac(\fX_k)$
surjects onto $\prod_{i=1}^r\Jac(\bar C_i)$ by Lemma~\ref{lem:rep.deg.g.avoid} below.
Therefore $\bar\phi_{E_2}$ is surjective for generic $\bar E_2$.  We fix 
a divisor $E_2$ whose reduction $\bar E_2$ is generic in this sense.

Now let $z\in J_0(K)$ and let $\bar z$ be its reduction in $\Jac(\fX_k)$.
By the above we know that there exists $E_1\in R\times C^{(*)}(K)$ such
that $\bar\phi_{E_2}(\bar E_1) = \bar z$.  Hence $\phi_{E_2}(E_1)$ and $z$
are both contained in $\red\inv(\bar z)$.  This formal fiber is isomorphic
to an open $g$-ball $\B(1)_+^g$; fix an isomorphism
$\red\inv(\bar z)\cong\B(1)_+^g$ such that $\phi_{E_2}(E_1)$ is identified
with zero.  Since $z$ is contained in some smaller closed ball
$\B(\rho)^g\subset\B(1)_+^g$, there exists a morphism
$\psi:\B(\rho)^g\to\red\inv(\bar z)$ such that 
$\psi(0) = \phi_{E_2}(E_1) = [E_1-E_2]$ and whose image contains $z$.  
By Bosch and L\"utkebohmert's Homotopy
Theorem~\cite[Theorem~3.5]{bosch_lutk:uniformizationII}, for generic 
$\bar E_2$ the map $\psi$ lifts through $R\times C^{(*)}$: that is, there
exists a map $\psi':\B(\rho)^g\to R\times C^{(*)}$ such that
$\phi_{E_2}\circ\psi' = \psi$.  Choosing $E_2$ so that $\bar E_2$
satisfies all necessary genericity conditions, we have that $z$ is
contained in the image of $\phi_{E_2}$, as desired.\qed

\smallskip
The following geometric lemma was used in the above proof.

\begin{lem} \label{lem:rep.deg.g.avoid}
  Let $k$ be an algebraically closed field and let $C$ be a smooth,
  proper, connected $k$-curve of genus $g$.  Let $D\subset C(k)$ be a
  finite set of points and let $\sL$ be a degree-zero line bundle on $C$.
  For a generic effective divisor $E_2\in\Div^g(C\setminus D)$, there
  exists an effective divisor $E_1\in\Div^g(C\setminus D)$ 
  such that $\sL\cong\sO(E_1-E_2)$.
\end{lem}

\pf Write $D = \{ x_1,\ldots, x_n\}$.  For any effective divisor
$E_2\in\Div^g(C)$ we have $h^0(\sL(E_2)) > 0$ by the Riemann-Roch
theorem, which means that $\sL$ admits a meromorphic section $s$ with
poles bounded by $E_2$; hence $\div(s) = E_1 - E_2$ for some effective
divisor $E_1\in\Div^g(C)$.  We need to show that for generic $E_2$, we can
choose $E_1$ with support disjoint from $D$.  It is enough to prove that
a generic effective divisor $E_2\in\Div^g(C\setminus D)$ satisfies
$h^0(\sL(E_2-(x_i)))=0$ for all $i=1,\ldots,n$ since then no global
section of $\sL(E_2)$ can have a zero at any $x_i$. 

A line bundle $\sL'$ on $C$ of degree $g-1$ is called \emph{special} if
$h^0(\sL')> 0$.  The special line bundles are exactly those contained in
the theta divisor $\Theta\subset\Pic^{g-1}(C)$, the image of the natural
map $C^{(g-1)}\to\Pic^{g-1}(C)$.  Hence we want to prove that for generic
$E_2$, the line bundle $\sL(E_2-(x_i))$ is not special for all
$i=1,\ldots,n$.  Let $\mu_i:\Pic^g(C)\isom\Pic^{g-1}(C)$ be the isomorphism
defined on $K$-points by $\sM\mapsto\sM(-x_i)$ and let
$U = \Pic^g(C)\setminus\bigcup_{i=1}^n \mu_i\inv(\Theta)$.  This is a dense
open subset of $\Pic^g(C)$, and for a line bundle $\sM\in U(K)$, by construction
$\sM(-x_i)$ is not special for all $i$.  Any divisor $E_2$ mapping to $U$
under the proper birational morphism $C^{(g)}\to\Pic^g(C)$ given by
$E\mapsto\sL(E)$ satisfies our requirements.\qed

\begin{rem} \label{rem:naive.false}
  With the notation in the above theorem,
  it is not true in general that there exist $D_1,D_2\in\Div^n(X)$ such
  that $D \coloneq D_1 - D_2$ is principal and $\tau_*(D) = \bar D$.
  For example, let $X$ be an elliptic curve with good reduction.  Then $X$
  has a skeleton $\Gamma$ which is a line segment, corresponding to a
  blow-up of its smooth model along a closed point on the special fiber.
  Let $P_1,P_2$ be the endpoints of $\Gamma$.  Then $(P_1)-(P_2)$ is a
  principal divisor on $\Gamma$, but a divisor of the form
  $(p_1)-(p_2)\in\Div^0(X)$ is principal if and only if $p_1=p_2$.
\end{rem}

\section{An application to tropical geometry}
\label{section:faithful_trop}

We conclude with an application of Theorem~\ref{thm:prin.surj} (or more precisely of Corollary~\ref{cor:merom.surj}) 
to ``faithful tropical representations'' of skeleta in the sense of \cite{bpr:trop_curves}.

\paragraph
In order to state the result, we briefly recall some terminology from tropical geometry.
Let $K$ be a complete and algebraically closed non-trivial non-Archimedean field.
Let $X$ be a variety over $K$ and let $Y$ be a toric variety over $K$ with dense torus $T = \Spec K[M]$, where $M \cong \ZZ^n$ is the character lattice of $T$.
If $f : X \dashrightarrow Y$ is any rational map which restricts to a
morphism $U \to T$ with $U$ a Zariski-dense open affine subset of $X$, we
can define a tropicalization map 
$\trop(f) : U^{\an} \to N_\R := \Hom(M,\RR) \cong \RR^n$ similarly
to~\parref{par:unif.basics.2} by setting  
\[
\left( \trop(f)(x) \right) (u) := -\log|\chi^u \circ f|_x,
\]
where $\chi^u$ is the character of $T$ corresponding to $u \in M$ and $| \cdot |_x$ is the multiplicative seminorm on $K[U]$ corresponding to $x \in U^{\an}$.
In particular, if $X$ is a smooth proper $K$-curve and $\Gamma$ is a skeleton of $X$, then $\Gamma$ is contained in $U^{\an}$ for any Zariski-dense open subset of $X$ and thus there is a well-defined restriction map 
\[
\trop(f)|_\Gamma ~:~ \Gamma \To N_\R.
\]

Since $\trop(f)(\Gamma)$ has rational slopes, there is a natural {\em lattice length metric} on $\trop(f)(\Gamma)$ (see \cite[Paragraph 6.2]{bpr:trop_curves}) and
since $\Gamma$ itself is naturally a metric space, it makes sense to ask whether $\trop(f)|_\Gamma$ is an {\em isometry}.  By \cite[Theorem 6.20]{bpr:trop_curves},
there exists a closed immersion $f : X \to Y_{\Delta}$ into some toric variety $Y_{\Delta}$ such that $\trop(f)|_\Gamma$ is an isometry.  However, the toric variety
$Y_{\Delta}$ constructed in {\em loc.~cit.}~has very large dimension in general.  Using Corollary~\ref{cor:merom.surj} we can cut down the dimension of $Y$ significantly,
at the expense of replacing the closed immersion $f$ with a rational map.

\begin{thm}
\label{thm:faithfultrop}
If $X$ is a smooth proper $K$-curve and $\Gamma$ is any skeleton of $X$, then there is a rational map $f : X \dashrightarrow {\mathbf P}^{3}$ such that the restriction of $\trop(f)$ to $\Gamma$ is an isometry onto its image (with respect to the lattice length metric on $\trop(f)(\Gamma) \subset {\mathbf R}^{3}$).
\end{thm}

\pf 
Let $G$ be a $\Lambda$-rational model for $\Gamma$ in the sense of~\parref{par:graph.basics}.
By Corollary~\ref{cor:merom.surj}, it suffices to prove that there are $\Lambda$-rational tropical meromorphic functions $F_1,F_2,F_3$ on $\Gamma$ which are linear on each edge of some $\Lambda$-rational subdivision $G''$ of $G$
such that (i) the $F_i$ separate points of $\Gamma$ and (ii) for each edge
$e$ of $G''$, the greatest common divisor of the absolute values of the slopes of the $F_i$ along $e$ is equal to~$1$.
(Condition (ii) guarantees that $\trop(f) : \Gamma \to {\mathbf R}^{3}$ is
a local isometry; c.f.~\cite[Definition~6.5 and~Remark 6.6]{bpr:trop_curves}, 
and condition (i) guarantees that it is injective.)

We define such functions $F_1,F_2,F_3$ as follows.  
Fix a loopless model $G$ for $\Gamma$ with no multiple edges and label the vertices of $G$ as $v_1,\ldots,v_n$ and the edges of $G$ as $e_1,\ldots,e_m$.  
Let $G'$ be the `barycentric subdivision' of $G$ in which a new vertex is inserted in the middle of each edge $e$, thereby replacing $e$ with two new edges of half
the original length.
For each edge $e=vw$ of $G$ and each $x \in e$, let $\hat{x}$ be the reflection of $x$ with respect to the midpoint $\tilde{e}$ of $e$.

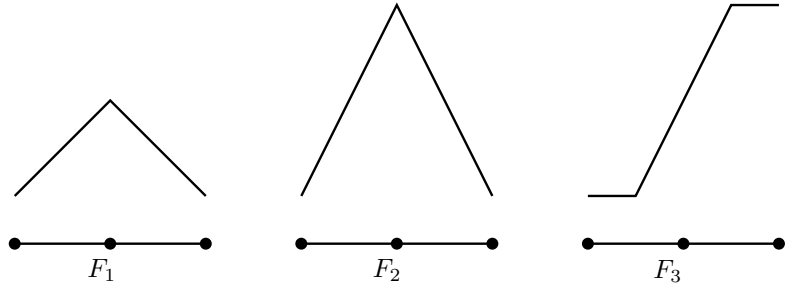
\genericfig[ht]{Fi}{Graphs of the functions $F_i$ in the proof of
  Theorem~\ref{thm:faithfultrop} along a fixed edge of $G$.} 

\begin{enumerate}
\item Let $F_1$ be the function which, for each edge $e=vw$ of $G$, increases linearly with slope $1$ from $v$ to $\tilde{e}$ and then decreases with slope $1$
from $\tilde{e}$ to $w$, with $F_1(v)=F_1(w)=0$ (so that $F_1$ is linear on each edge of $G'$).
\item Choose distinct positive integers $\sigma_1,\ldots,\sigma_m$ and let $F_2$ be the function which, for each edge $e_i=vw$ of $G$, increases linearly with slope 
$\sigma_i$ from $v$ to $\tilde{e_i}$ and then decreases with slope $\sigma_i$
from $\tilde{e_i}$ to $w$, with $F_1(v)=F_1(w)=0$ (so that $F_2$ is linear on each edge of $G'$).
\item Choose distinct positive real numbers $\alpha_1, \ldots, \alpha_n \in \Lambda$ and let $F_3$ be any $\Lambda$-rational tropical meromorphic function on $\Gamma$ 
such that $F_3(v_i)=\alpha_i$ for all $i$, and which is {\em asymmetric} in the sense that $F_3(x) \neq F_3(\hat{x})$ for all $x \in \Gamma \backslash V(G')$.
(For example, one can choose $F_3$ so that along each edge $e=vw$ of $G$,
$F_3$ has slope zero in some neighborhoods of $v$ and $w$ and slope $\pm \mu$ in some 
neighborhood of $\tilde{e}$ for an integer $\mu$ with $|\mu|$ sufficiently large.)
\end{enumerate}

See Figure~\ref{fig:Fi}.
Let $G''$ be a refinement of $G'$ such that $F_3$ is linear along each edge of $G''$.  We claim that $F_1,F_2,F_3$ satisfy conditions (i) and (ii) above.
Condition (ii) is clear, since $F_1$ has slope $\pm 1$ along each edge of $G''$.
To verify condition (i), suppose $x,y \in \Gamma$ and that $F_i(x)=F_i(y)$ for $i=1,2,3$.  We want to show that $x=y$.
Suppose $x \in e_i$ and $y \in e_j$, and identify $e_i$ (resp. $e_j$) with $[0,\ell(e_i)]$ (resp. $[0,\ell(e_j)]$) via an isometry $\iota_i$ (resp. $\iota_j$) 
in such a way that $\iota_i(x) \in [0,\ell(e_i)/2]$ (resp. $\iota_j(y) \in [0,\ell(e_j)/2]$).
Since $F_1(x)=F_1(y)$, we must have $\iota_i(x) = \iota_j(y)$.
Since $F_2(x)=F_2(y)$, we must have $\sigma_i \iota_i(x) = \sigma_j \iota_j(y)$, and therefore either $i=j$ or $\iota_i(x)=\iota_j(y)=0$ (i.e., $x,y$ are vertices of $G$).
In the latter case, since $F_3$ takes distinct values on distinct vertices of $G$ we conclude that $x=y$ as desired.
So suppose that $i=j$ and that $x,y$ are not vertices of $G$.  From the relation $\iota_i(x) = \iota_j(y)$ we conclude that either $x=y$ or $x=\hat{y}$.
Since $F_3$ is asymmetric, the relation $F_3(x)=F_3(y)$ implies that $x=y=\tilde{e_i}$.
\qed

\begin{rem}
The proof of Theorem~\ref{thm:faithfultrop} implies  the purely combinatorial fact that any metric graph $\Gamma$ admits a quasi-balanced isometric embedding in $\RR^3$, where `quasi-balanced' means that one can add infinite rays (with multiplicities) to $\Gamma$ in such a way that the resulting one-dimensional weighted polyhedral
complex $\tilde{\Gamma}$ satisfies the balancing condition from tropical
geometry.  (See for example~\cite[\S3.4]{maclagan_sturmfels:book} for a
discussion of the balancing condition.)
\end{rem}

\begin{rem}
If $X$ is an elliptic curve with bad reduction and $\Gamma$ is the minimal skeleton of $X$ (which is homeomorphic to a circle), 
then by~\cite[Theorem~7.2]{bpr:trop_curves}
or~\cite{chan_sturmfels:honeycomb} there is a closed immersion $f : X \to \PP^2$ such that the restriction of $\trop(f)$ to $\Gamma$ is an isometry.
However, it is clear that such a result cannot hold for arbitrary curves $X$ since (a) not every algebraic curve can be embedded in $\PP^2$ and (b) not every
graph is planar.
\end{rem}

\begin{rem}
We conjecture that Theorem~\ref{thm:faithfultrop} can be strengthened by requiring $f$ to be a {\em closed immersion} of $X$ in a 3-dimensional projective toric variety 
$Y_{\Delta}$ (in which case $\Gamma$ is ``faithfully represented'' in the sense of \cite[Paragraph 6.15.2]{bpr:trop_curves}).  One might even be able to always take 
$Y_{\Delta} = \PP^3$.
\end{rem}

\begin{rem}
Corollary~\ref{cor:merom.surj} can also be used to streamline the proof of \cite[Theorem 6.20]{bpr:trop_curves}, reducing the number of rational functions required
to get a faithful representation.  In particular, Corollary~\ref{cor:merom.surj} obviates the need for Lemmas 6.18 and 6.19 in \cite{bpr:trop_curves}, and provides a 
powerful generalization of such results.
\end{rem}

\bibliographystyle{thesis}
\bibliography{harm_etale}

\end{document}